\documentclass[12pt]{article}
\usepackage{amsmath}
\usepackage{amssymb}
\usepackage[ctr1,eng]{myarticl}

\newcommand{\alp }{\alpha }
\newcommand{\bet }{\beta }
\newcommand{\antip }{\kappa }
\newcommand{\cA }{\mathcal{A}}
\newcommand{\cB }{\mathcal{B}}
\newcommand{\cF }{\mathcal{F}}

\newcommand{\cH }{\mathcal{H}}
\newcommand{\cI }{\mathcal{I}}

\newcommand{\copr }{\varDelta }
\newcommand{\End }{\mathrm{End}}

\newcommand{\hght }[1]{h_{#1}}
\newcommand{\lact }{.}
\newcommand{\lch }[1]{#1{}_{\scriptscriptstyle \mathrm{L}}}
\newcommand{\lchl }{\ell _{\scriptscriptstyle \mathrm{L}}}
\newcommand{\lcoa }{\delta }
\newcommand{\lgf }[1]{#1{}^{\scriptscriptstyle \mathrm{L}}}
\newcommand{\lgfl }{\ell ^{\scriptscriptstyle \mathrm{L}}}
\newcommand{\lgh }{`L\mbox{'}}
\newcommand{\mbbX }{\mathbb{X}}
\newcommand{\Nbar }{\bar{N}}
\newcommand{\Ndbasis }{\boldsymbol{\mathrm{e}}}
\newcommand{\ndN }{\mathbb{N}}
\newcommand{\ndQ }{\mathbb{Q}}
\newcommand{\ndZ }{\mathbb{Z}}
\newcommand{\ot }{\otimes }

\newcommand{\paar }[2]{\langle #1,#2\rangle }
\newcommand{\paarb }[2]{\bigg\langle #1,#2\bigg\rangle }

\newcommand{\qfact }[2]{(#1)^!_{#2}}
\newcommand{\qnum }[2]{(#1)_{#2}}
\newcommand{\rch }[1]{#1{}_{\scriptscriptstyle \mathrm{R}}}
\newcommand{\rchl }{\ell _{\scriptscriptstyle \mathrm{R}}}

\newcommand{\rgf }[1]{#1{}^{\scriptscriptstyle \mathrm{R}}}
\newcommand{\rgfl }{\ell ^{\scriptscriptstyle \mathrm{R}}}
\newcommand{\rgh }{`R\mbox{'}}
\newcommand{\roo }{\mathrm{\mathbf{r}}}
\newcommand{\spaar }[2]{\langle \iota (\tau (#1)),\tau (#2)\rangle }
\newcommand{\StBr }{\sigma }

\newcommand{\suplet}[1]{[#1]}
\newcommand{\YD }{Yetter--Drinfel'd }
\newcommand{\YDcat }{\mathcal{YD}}

\title{Finite dimensional rank 2 Nichols algebras\\
of diagonal type I: Examples}
\author{I.~Heckenberger\thanks{email:
Istvan.Heckenberger@math.uni-leipzig.de}}

\begin{document}
\maketitle

\begin{abstract}
Nichols algebras naturally appear in the classification
of finite dimensional pointed Hopf algebras. Assuming only that the base
field has characteristic zero several new finite dimensional rank
2 Nichols algebras of
diagonal type are listed. Each of them is described in terms of generators and
relations. A Poincar\'e--Birkhoff--Witt basis of all of these Nichols algebras
is given and their dimension is computed. Together with the continuation
of this paper this gives an answer to a question of Andruskiewitsch.

Key Words: Hopf algebra, Lyndon words, full binary tree

MSC2000: 17B37, 16W35
\end{abstract}

\section{Introduction}

Nichols algebras are currently studied intensively as part of Hopf
algebra theory, see e.\,g.~\cite{a-Nichols78}, \cite{a-AndrSchn98},
\cite{a-AndrGr99}, \cite{inp-Grana99}, \cite{a-Khar99}, \cite{inp-AndrSchn02},
\cite{a-Ufer03}, \cite{inp-Andr03} and the references therein.
They can also be seen as a generalization of exterior
algebras \cite{a-Woro2}.
One of the most remarkable property of Nichols algebras is their
relation to root systems of semisimple Lie algebras \cite{a-AndrSchn00}.

The classification of finite dimensional Hopf algebras is one of the oldest
problems of Hopf algebra theory. As part of it also pointed Hopf algebras
are studied. In their pioneering work
\cite{a-AndrSchn98} N.~Andruskiewitsch and H.-J.~Schneider propose an
elegant method to carry out a classification for this special class of Hopf
algebras. The method consists of several steps. It relies on the fact
that to any pointed Hopf algebra $\cH $ one can associate a graded Hopf
algebra $\mathrm{gr}\,\cH $ corresponding to the coradical filtration of
$\cH $. Further, it is known that $\mathrm{gr}\,\cH $ is isomorphic to the
Radford biproduct $R\#kG$ of a braided Hopf algebra $R$ and the group
algebra $kG$ (which is the coradical of $\cH $) where $G$ is the group
of group-like elements of $\cH $.
Thus first one has to determine all possible pairs $(R,G)$ and then
one has to construct all Hopf algebras (the so called liftings) such that
the corresponding graded Hopf algebra is $R\#kG$. Usually $R$ is a Nichols
algebra, i.\,e.~it is generated by the vector space $V$ of its (twisted)
primitive elements. The space $V$ turns out to be a \YD module which
determines the Nichols algebra $R$ (which is then denoted by $\cB (V)$)
uniquely. There exists an explicit description of $\cB (V)$
by Schauenburg \cite{a-Schauen96} in terms of $V$ and $G$. Nichols algebras
appeared in this form also as quantized exterior algebras of Woronowicz
in \cite{a-Woro2}.

If the action and coaction of $kG$ on $V$ are simultaneously diagonalizable
then $V$ is said to be of diagonal type.
Currently there exist
several classification results in this case which tell that under some
assumptions on $G$ all \YD modules are related to symmetrizable Cartan
matrices of finite type. Further there exists a list of examples which do not
fit into the above classification scheme \cite{inp-AndrSchn02}.
Andruskiewitsch stated in \cite{inp-Andr02} the following question.

\addvspace{\baselineskip}

\textsc{Question} 5.40. Given a braided vector space $V$ of diagonal type
and dimension 2, decide when $\cB (V)$ is finite dimensional. If so,
compute $\dim \cB (V)$, and give a ``nice'' presentation by generators
and relations.

\addvspace{\baselineskip}

The aim of the present paper is to give a complete answer to
Question 5.40 of Andruskiewitsch.
In the first part examples of rank 2 Nichols algebras $\cB (V)$
are listed. All of them are related to one of 22 different full
binary trees. {}From these trees one can also read off a
generating set of relations and a Poincar\'e--Birkhoff--Witt basis
of $\cB (V)$. This allows the computation of the dimension
of $\cB (V)$. In the continuation of this paper the first part of
Question 5.40 will be addressed.

Our main goal is to describe rank 2 Nichols algebras with help of
full binary trees. The latter can be considered as analogues of
the root system of a rank 2 semisimple Lie algebra. Indeed, if the
braiding of $V$ is of Cartan type, then the set $N_2(T)$ of those nodes
of the full binary tree $T$ which have two children corresponds to
the set of nonsimple positive roots of the Lie algebra. In general,
the edges of $T$ and the set of nodes having no children correspond
to relations of the Nichols algebra.

The list of (all known and) new examples is contained in Theorem \ref{t-class}.
Their construction uses several ideas. The computational
part needs heavily the fact that there exists an action of $\cB (V^*)\#kG$
on $\cB (V)$, see Lemma \ref{l-nichpaarung} and Corollary \ref{f-nichpaarung}.
This already seems to be known and there exist various forms
of it in the literature, usually as some bilinear pairing \cite{b-Lusztig93}
or as quantum differential operators \cite{a-Grana00}.
The theoretical part is based on the one side on an old result of
Stern \cite{a-Stern1858} on some special sequences of pairs of integer
numbers. This theory is part of graph theory and is contained also in the
modern literature \cite{b-CLR90}. An adapted version of it is described 
in Section \ref{ss-rank2types}. On the other hand deep results of Kharchenko
\cite{a-Khar99} on the structure of general Nichols algebras are used.
They say that any Nichols algebra of rank $n$ has a
Poincar\'e--Birkhoff--Witt basis which corresponds
to a subset of all Lyndon words of an alphabet with $n$ letters.  
Further Kharchenko also proves very strong restrictions on the generating set
of relations of the Nichols algebra.
Finally, using the special situation when $V$ has rank two one can 
relate full binary trees $T$ and Nichols algebras $\cB (V)$ such that
nodes of $T$ correspond to PBW generators and relations of $\cB (V)$.
This is done in Section \ref{ss-Lyndonfull}.
In order to check correlations between $T$ and $V$ one still has to perform
tedious computations but in advantage one can
eventually avoid the use of computer algebra programs.

If not stated otherwise the definitions and notation follow
\cite{a-AndrSchn98}.
Throughout this paper $k$ denotes a field of characteristic zero and
tensor products $\ot $ are taken over this field. For Hopf algebras
the coproduct and the antipode are denoted by $\copr $
and $\kappa $, respectively. We use the Sweedler notation
$\copr (a)=a_{(1)}\ot a_{(2)}$ for elements $a$ of a Hopf algebra. The set of
natural numbers not including 0 is denoted by $\ndN $ and we write
$\ndN _0=\ndN \cup \{0\}$.

The author wants to thank A.~Joseph and S.~Ufer for stimulating discussions
and N.~Andruskiewitsch and M.~Gra{\~n}a for helpful remarks.

\section{Nichols algebras}
\label{s-Nichols}

\subsection{The dual of a Nichols algebra}
\label{ss-pairing}

Suppose that $k$ is a field of characteristic zero,
$G$ an abelian group, and $V\in {}^{kG}_{kG}\YDcat $ a finite dimensional
\YD module with completely reducible $kG$-action. Let $\lcoa :V\to kG\ot V$
and $\lact :kG \ot V\to V$ denote the left coaction and left action of $kG$
on $V$, respectively. If $G$ is finite and $k$ is algebraically closed then
the condition on complete reducibility is automatically fulfilled. Anyway,
in such a case the braiding $\sigma \in \End _k(V\ot V)$ of $V$ where
\begin{align*}
\sigma (v\ot w)=&(v_{(-1)}\lact w)\ot v_{(0)},&
\sigma ^{-1}(v\ot w)=&w_{(0)}\ot (\antip ^{-1}(w_{(-1)})\lact v),
\end{align*}
and $\lcoa (v)=v_{(-1)}\ot v_{(0)}$ for $v\in V$, is called \textit{of diagonal
type}. Let $\cB (V)$ denote the Nichols algebra generated by $V$.
More precisely, as proved in \cite{a-Schauen96} and noted in
\cite[Prop.~2.11]{inp-AndrSchn02},
\begin{align*}
\cB (V)=k\oplus V\oplus \bigoplus _{m=2}^\infty V^{\ot m}/\ker S_m
\end{align*}
where $S_m\in \End _k(V^{\ot m})$, $S_{1,j}\in \End _k(V^{\ot j+1})$,
\begin{align*}
S_m&=\prod _{j=1}^{m-1}(\id ^{\ot m-j-1}\ot S_{1,j}),\\
S_{1,j}&=\id +\sigma ^{-1}_{12}+\sigma ^{-1}_{12}\sigma ^{-1}_{23}
+\cdots +\sigma ^{-1}_{12}\sigma ^{-1}_{23}\cdots \sigma ^{-1}_{j,j+1}
\end{align*}
(in leg notation) for $m\ge 2$ and $j\in \ndN _0$.
Let $\cB (V)^+$ denote the unique maximal ideal of $\cB (V)$.

Now we are going to give a slightly modified version
(see Lemma \ref{l-nichpaarung}) of Lusztig's bilinear form \cite{b-Lusztig93},
cf.~\cite{inp-AndrSchn02}. It is closely related to
differential operators on $\cB (V)$ defined in \cite{inp-Grana99},
\cite{a-Grana00} and \cite{inp-AndrSchn02}.

Let $V^*$ denote the \YD module dual to $V$. More precisely, for $f\in V^*$
one has
\begin{align*}
(h\lact f)(v)&=f(\antip (h)\lact v),& f_{(0)}(v)f_{(-1)}&=f(v_{(0)})
\antip ^{-1}(v_{(-1)}) && \forall v\in V, h\in kG,
\end{align*}
where $\lcoa (f)=f_{(-1)}\ot f_{(0)}$. Note that there exists a linear map
$\paar{\cdot }{\cdot } :V^*\times \cB (V)\to \cB (V)$ such that
$\paar{f}{1}:=0$ and
\begin{align*}
\paar{f}{\rho }&:=\sum _if(a_i)\rho _i,\text{ where $S_{1,m-1}(\rho )
=\sum _ia_i\ot \rho _i\in V\ot V^{\ot {m-1}}$}
\end{align*}
for $m>0$ and $\rho \in V^{\ot m}$. Since $(f\ot \id )\sigma ^{-1}(v\ot w)=
(f_{(-1)}\lact v)f_{(0)}(w)$ for all $v,w\in V$ and $f\in V^*$ the map
$\paar{\cdot }{\cdot }$ has the property
\begin{align*}
\paar{f}{\rho \rho '}=\paar{f}{\rho }\rho '+(f_{(-1)}\lact \rho )
\paar{f_{(0)}}{\rho '}
\end{align*}
for all $f\in V^*$ and $\rho ,\rho '\in \cB (V)$.
Note that $\rho =0$ in $\cB (V)$ if and only if $\paar{f}{\rho }=0$ for all
$f\in V^*$. Moreover, the mappings $D_i$ in \cite{inp-AndrSchn02} and
\cite{inp-Grana99} can be considered as mappings $\paar{y_i}{\cdot }$ where
$\{y_i\}$ is a canonical basis of $V^*$ (see the notation in Subsection
\ref{ss-conventions}).

Let $\cB (V^*)\# kG$ denote the set $\cB (V^*)\ot kG$ with the product
\begin{align*}
(f'\ot g')(f''\ot g''):=f'(g'\lact f'')\ot g'g''
\end{align*}
for all $f',f''\in \cB (V^*)$ and $g',g''\in G$. Then $\cB (V^*)\# kG$
becomes a Hopf algebra with coproduct
\begin{align*}
\copr (f)&=f\ot 1+\lcoa (f),& \copr (g)=g\ot g,& &f\in V^*,g\in G.
\end{align*}

\begin{lemma}\label{l-nichpaarung}
There exists a unique bilinear map
$\paar {\cdot }{\cdot }: (\cB (V^*)\# kG)\times \cB (V)\to \cB (V)$
satisfying $\paar{f}{v}=f(v)$ for $f\in V^*$, $v\in V$, and
\begin{gather*}
\paar{f}{\rho \rho '}=\paar{f}{\rho }\rho '+(f_{(-1)}\lact \rho )
\paar{f_{(0)}}{\rho '},\qquad
\paar{g}{\rho }=g\lact \rho ,\\
\paar{h_1h_2}{\rho }=\paar{h_1}{\paar{h_2}{\rho }}
\end{gather*}
for all $f\in V^*$, $g\in kG $, $h_1,h_2\in \cB (V^*)\# kG$, and
$\rho ,\rho '\in \cB (V)$.
\end{lemma}

Since $\copr (f)=f\ot 1+\lcoa (f)$ for $f\in V$ and $\copr (g)=g\ot g$
for $g\in G$ the following corollary is an immediate consequence of Lemma
\ref{l-nichpaarung}.

\begin{folg}\label{f-nichpaarung}
For any $f\in \cB (V^*)\# kG$ and $\rho ,\rho '\in \cB (V)$ the formula
$\paar{f}{\rho \rho '}=\paar{f_{(1)}}{\rho }\paar{f_{(2)}}{\rho '}$ holds
where $\copr (f)=f_{(1)}\ot f_{(2)}$.
\end{folg}

\begin{bew}[ of the Lemma]
It is clear that the restriction of $\paar{\cdot }{\cdot }$ to $V^*\times
\cB (V)$ has to be the map $\paar{\cdot }{\cdot }$ introduced above. Therefore
the uniqueness assertion immediately follows. Further, equation
\begin{align*}
\paar{g\lact f}{\rho }=g^{-1}\lact \paar{f}{g\lact \rho },\quad
g\in G,f\in V^*,\rho \in \cB (V),
\end{align*}
implies that there exists at least a map $\paar{\cdot }{\cdot }:
\left(\bigoplus _{m=0}^\infty V^{\ot m}\# kG\right)\times \cB (V)\to
\cB (V)$ satisfying the above equations.

For $1\le l\le m$ there exist linear maps $S_{l,m-l}\in \End _k(V^{\ot m})$
such that $\prod _{j=1}^l(\id ^{\ot l-j}\ot S_{1,m-l+j-1})=
(S_l\ot \id ^{\ot m-l})S_{l,m-l}$.
Thus for $l,m\in \ndN _0$, $l\le m$, $f_i\in V^*$,
$1\le i\le l$, and $\rho \in V^{\ot m}$ one has
\begin{align*}
\paar{f_1}{\paar{f_2}{\cdots \paar{f_l}{\rho }}}=
(f_l\ot f_{l-1}\ot \cdots f_1\ot \id ^{m-l})((S_l\ot \id ^{m-l})
S_{l,m-l}(\rho )).
\end{align*}
Moreover, direct computation shows that
$(f\ot g)(\sigma ^{-1}(v\ot w))=\sigma ^{-1}(g\ot f)(w\ot v)$ for $f,g\in V^*$
and $v,w\in V$.
Therefore $\paar{\ker S_l^*}{\cB (V)}=0$ where
$S_l^*\in \End _k(V^*{}^{\ot l})$ is defined analogously to $S_l$.
\end{bew}

\subsection{Conventions}
\label{ss-conventions}

Let $d\in \ndN $, $g_i\in kG$ for $1\le i\le d$, $q_{ij}\in k$
for $i,j\in \{1,2,\ldots ,d\}$, and
$\{x_i\,|\,1\le i\le d\}$ a basis of $V$ such that
\begin{align*}
\lcoa (x_i)&=g_i\ot x_i,&  g_i\lact x_j&=q_{ij}x_j.
\end{align*}
Such a basis always exists and is called \textit{a canonical basis} of $V$.
Let $\{y_i\,|\,1\le i\le d\}$ denote the dual basis of $V^*$.
Then
\begin{align*}
\lcoa (y_i)&=g_i^{-1}\ot y_i,& g_i\lact y_j&=q_{ij}^{-1}y_j,&
\sigma (y_i\ot y_j)&=q_{ij}y_j\ot y_i.
\end{align*}
Thus for diagonal braidings the linear map $\iota :V\to V^*$,
$\iota (x_i):=y_i$ for $1\le i\le d$, extends to an algebra isomorphism
$\iota :\cB (V)\to \cB (V^*)$.

Let $\{\Ndbasis _i\,|\,1\le i\le d\}$ denote a basis of the $\ndN _0$-module
$\ndN _0^d$. The tensor algebra $V^\ot =\bigoplus _{m=0}^\infty V^{\ot m}$
admits an $\ndN _0^d$-grading defined by $\deg x_i:=\Ndbasis _i$,
$1\le i\le d$. The corresponding total degree $\mathrm{totdeg}$ is the
$\ndN _0$-grading of $V^\ot $ defined by $\mathrm{totdeg}(x_i)=1$ for
$1\le i\le d$. Note that the linear maps
$\paar{y_i}{\cdot }$, $1\le i\le d$, are $\ndN _0^d$-graded and hence
the Nichols algebra $\cB (V)$ inherits the $\ndN _0^d$- and $\ndN _0$-grading
of $V^\ot $.
Let $\cB (V)_n$, $n\in \ndN _0$, denote the set of homogeneous elements
of $\cB (V)$ of total degree $n$.

To a given \YD module $V$ as in this subsection there exist a unique group
homomorphism $g:\ndZ ^d\to G$ and a unique bicharacter
$\chi :\ndZ ^d\times \ndZ ^d\to k$ satisfying
$g(\Ndbasis _i)=g_i$ and $\chi (\Ndbasis _i,\Ndbasis _j)=q_{ij}$ for
$i,j\in \{1,2,\ldots ,d\}$.
For notational convenience we will also write $g(x)$ and $\chi (x',x'')$
instead of $g(\deg x)$ and $\chi (\deg x',\deg x'')$ for homogeneous
elements $x,x',x''\in V^\ot $ and $x,x',x''\in \cB (V)$.
Note that if there exist $i,j$ such that
$q_{ij}\not=q_{ji}$ then the bicharacter $\chi $ is not symmetric!

\section{Rank 2 Nichols algebras}
\label{s-rank2}

\subsection{Types of Nichols algebras}
\label{ss-rank2types}

For basic definitions in this section we refer to \cite{b-CLR90} and
\cite{b-Loth83}.
Recall that a binary tree $T$ is a (nonempty finite) tree such that each node
has at most two children. One says that $T$ is \textit{full} if each node of
$T$ has exactly zero or two children \cite{b-CLR90}. For examples see Appendix
\ref{app-types}. For the set of nodes of a full binary tree $T$ which have
zero and two children, respectively, we use the symbol $N_0(T)$ and $N_2(T)$,
respectively. Let $\roo (T)$ or simply $\roo $ denote the root of the binary
tree $T$. Further, we write $N(T)=N_0(T)\cup N_2(T)$ for the set of all
nodes of $T$. Let $\{\lgh ,\rgh \}$ be a set with two elements and define
$\Nbar _2(T):=N_2(T)\cup \{\lgh ,\rgh \}$,
$\Nbar (T):=N(T)\cup \{\lgh ,\rgh \}$ (disjoint unions).

\begin{defin}
Let $T$ be a full binary tree and $a\in N(T)$.
A node $b\in \Nbar _2(T)$ is called the \textit{left godfather of $a$},
denoted by $b:=\lgf{a}$, if one of the following conditions holds.
\begin{itemize}
\item
$a=\roo (T)$ and $b=\lgh $,
\item
$a$ is the right child of $b$,
\item
$a$ is the left child of its parent $c$ and $b$ is the left godfather of $c$.
\end{itemize}
Similarly one defines the \textit{right godfather $\rgf{a}$ of $a$} by
replacing everywhere left by right and vice versa and setting
$\rgf{\roo }:=\rgh $. If $a\in N_2(T)$ then let
$\lch{a}$ and $\rch{a}$ denote the left and right child of $a$, respectively.
\end{defin}
Note that $\lgf{\cdot }$ and $\rgf{\cdot }$ are well defined maps from $N(T)$
to $\Nbar _2(T)$ and any $a\in N(T)$ is uniquely determined by $\lgf{a}$ and
$\rgf{a}$. Indeed, $a\in N(T)$, $b=\lgf{a}$ implies that there is a subset
$\{a_1,a_2,\ldots ,a_m\}$ of $N(T)$ such that $m\in \ndN $, $a_m=a$,
$a_{i+1}$ is the left child of $a_i$ for all
$i>0$, and either $a_1=\rch{b}$ or $a_1=\roo $, $b=\lgh $.
Further, $\rgf{a_{i+1}}=a_i$ for all $i>1$ and
$\rgf{a_1}\notin \{a_1,a_2,\ldots ,a_m\}$. Thus $a$ is uniquely determined by
$\lgf{a}$ and $\rgf{a}$.

We define functions $\lgfl ,\rgfl ,\lchl $, and $\rchl :N(T)\to \ndN $
recursively to denote lengths of certain branches of $T$. Set
\begin{align*}
\lgfl (a):=&\begin{cases}
\lgfl (\lgf{a})+1 & \text{if $\rch{\lgf{a}}=a$,}\\
1 & \text{else,}
\end{cases}
&
\lchl (a):=&\begin{cases}
1 & \text{if $a\in N_0(T)$,}\\
\lchl (\lch{a})+1 & \text{else,}
\end{cases}
\\
\rgfl (a):=&\begin{cases}
\rgfl (\rgf{a})+1 & \text{if $\lch{\rgf{a}}=a$,}\\
1 & \text{else,}
\end{cases}
&
\rchl (a):=&\begin{cases}
1 & \text{if $a\in N_0(T)$,}\\
\rchl (\rch{a})+1 & \text{else.}
\end{cases}
\end{align*}

Any full binary tree $T$ can be identified with a subtree of the infinite
\textit{Stern--Brocot tree} \cite[pp.~116--117]{b-GKP94}, see also
\cite{a-Stern1858}.
This means that there exists a map $\StBr _T:\Nbar (T)\to \ndZ \times \ndZ $
defined recursively by $\StBr _T(\lgh ):=(0,1)$, $\StBr _T(\rgh ):=(1,0)$, and
$\StBr _T(a):=\StBr _T(\lgf{a})+\StBr _T(\rgf{a})$ for any $a\in N(T)$.
Note that $\StBr _T(a)\in \ndN \times \ndN $ for $a\in N(T)$ by the definition
of $\lgf{\cdot }$ and $\rgf{\cdot }$ and since $\StBr _T(\roo )=(1,1)$.
Thus the map $Q:\Nbar (T)\to \ndQ \cup \{+\infty \}$,
\begin{align*}
Q(a):=\begin{cases}
r/s & \text{if $\StBr _T(a)=(r,s)$, $s\not=0$,}\\
+\,\infty & \text{if $a=\rgh $,}
\end{cases}
\end{align*}
and the total order $<$ on $\ndQ $ induce an order $<_Q$
on $\Nbar (T)$ such that for all $a\in N(T)$ the relations 
$\lgh <_Qa<_Q\rgh $ hold.
There is another natural map $|\StBr _T|:\Nbar (T)\to \ndZ $
defined by $|\StBr _T|(a)=r+s$ whenever $\StBr _T(a)=(r,s)$. It will be used
mainly for inductive proofs.

Assertions (i)--(iii) of the following Lemma were proved e.\,g.\ in
\cite{b-GKP94}.

\begin{lemma}\label{l-fbtreeord}
Let $T$ be a full binary tree and $a,b\in N(T)$.\\
(i) $\StBr _T(a)=(r,s)$, $\StBr _T(\lgf{a})=(r_1,s_1)$,
$\StBr _T(\rgf{a})=(r_2,s_2)$ $\Rightarrow $ $rs_1-r_1s=r_2s-rs_2
=r_2s_1-r_1s_2=1$.\\
(ii) The entries of $\StBr _T(c)$, $c\in \Nbar (T)$, are relatively prime
integers.\\
(iii) The map $Q:\Nbar (T)\to \ndQ \cup \{+\infty \}$ is injective and hence
$<_Q$ is a total order on $\Nbar (T)$.\\
(iv) $\lgf{a}<_Qa<_Q\rgf{a}$.\\
(v) $a<_Qb$, $|\StBr _T|(a)<|\StBr _T|(b)+|\StBr _T|(\lgf{b})$ $\Rightarrow $
$a\le _Q\lgf{b}$.\\
(vi) $a<_Qb$, $|\StBr _T|(b)<|\StBr _T|(a)+|\StBr _T|(\rgf{a})$
$\Rightarrow $ $\rgf{a}\le _Q b$.\\
(vii) $\lgf{a}\in N(T)$ $\Rightarrow $ $a<_Q\rgf{\lgf{a}}$.\\
(viii) $\rgf{a}\in N(T)$ $\Rightarrow $ $\lgf{\rgf{a}}<_Qa$.\\
(ix) $a\in N(T), b\in \Nbar (T)$, $|\StBr _T|(b)\le |\StBr _T|(a)$,
$\lgf{a}<_Qb<_Q\rgf{a}$ $\Rightarrow $ $b=a$.
(x) If $a,b\in \Nbar _2(T)$ such that $a<_Qb$ then there exists
$c\in N(T)$ such that $a<_Qc<_Qb$.
\end{lemma}

\begin{bew}
(i) We prove this by induction on $|\StBr _T|(a)$. If $a=\roo $ then
$\StBr _T(a)=(1,1)$ and the assertion holds. Otherwise let $c\in N(T)$ denote the
parent of $a$. If $a$ is the left child of $c$ then $\rgf{a}=c$,
$\lgf{a}=\lgf{c}$, and $\StBr _T(c)=(r_2,s_2)$. Since
$|\StBr _T|(c)<|\StBr _T|(a)$ equation $\StBr _T(\lgf{c})=(r_1,s_1)$ and
induction hypothesis give
$r_2s_1-r_1s_2=1$. Thus $rs_1-r_1s=(r_1+r_2)s_1-r_1(s_1+s_2)=1$ and
$r_2s-rs_2=r_2(s_1+s_2)-(r_1+r_2)s_2=r_2s_1-r_1s_2=1$.
If $a=\rch{c}$ then one argues similarly.

(ii) This one gets from (i) using $\StBr _T(\lgh )=(0,1)$,
$\StBr _T(\rgh )=(1,0)$.

(iv) This follows from (i) and the fact that $\StBr _T(c)\in \ndN \times \ndN $
for $c\in N(T)$.

(v) Suppose to the contrary that $\lgf{b}<_Q a<_Q b$ and $|\StBr _T|(a)<|\StBr _T|(b)
+|\StBr _T|(\lgf{b})$. If $\StBr _T(a)=(r,s)$, $\StBr _T(b)=(r_2,s_2)$, and
$\StBr _T(\lgf{b})=(r_1,s_1)$ then $r_2s_1-r_1s_2=1$, $rs_1-r_1s\ge 1$
and $r_2s-rs_2\ge 1$. Therefore
\begin{align*}
(r_2s)s_1\ge (rs_1)s_2+s_1\ge r_1ss_2+s_2+s_1 \Rightarrow s_1+s_2\le s,\\
r_2(rs_1)\ge r_1(r_2s)+r_2\ge r_1rs_2+r_1+r_2 \Rightarrow r_1+r_2\le r.
\end{align*}
This is a contradiction to
$r+s=|\StBr _T|(a)<|\StBr _T|(b)+|\StBr _T|(\lgf{b})
=r_1+r_2+s_1+s_2$.

(vi) Use arguments as in (v).

(iii) If $Q(a)=Q(b)$ then $\StBr _T(a)=\StBr _T(b)$ by (ii).
Now $\lgf{b}<_Qa$, $|\StBr _T|(\lgf{b})<|\StBr _T|(b)=|\StBr _T|(a)$,
and (v) imply that $\lgf{b}\le _Q\lgf{a}$. By symmetry one gets
$Q(\lgf{a})=Q(\lgf{b})$ and similarly $Q(\rgf{a})=Q(\rgf{b})$.
Thus using that $\lgf{a}$ and $\rgf{a}$ determine $a$ uniquely
induction on $|\StBr _T|(a)$ gives the assertion.

(vii) Again we use induction on $|\StBr _T|(a)$.
Note that $a\not=\roo $ since $\lgf{a}\in N(T)$.
If $a$ is the left child of its parent $c$ then $\lgf{a}=\lgf{c}$
and $\rgf{a}=c$. By induction hypothesis $c<_Q\rgf{\lgf{c}}$ and hence (iv)
implies that $a<_Q\rgf{a}=c<_Q\rgf{\lgf{a}}$.
If on the other hand $a$ is the right child of its parent $c$ then
one gets $\lgf{a}=c$ and $\rgf{a}=\rgf{c}$. Thus by (iv) one obtains
$a<_Q\rgf{a}=\rgf{c}=\rgf{\lgf{a}}$.

(viii) The proof goes as for (vii).

(ix) Relation $b<_Q\rgf{a}$ and the converse of (vi) imply that $b\le _Qa$.
Similarly, $\lgf{a}<_Qb$ and the converse of (v) yield that $a\le _Qb$.
By (iii) one gets $a=b$.

(x) Set $c_1:=\rch{a}$ and $c_2:=\lch{b}$.
Then $a<_Qc_1$ and $c_2<_Qb$ by (iv). If $c_1=b$ then
$a=\lgf{b}=\lgf{c_2}<_Qc_2<_Q\rgf{c_2}=b$ again by (iv). In this case one can
take $c=c_2$.
Similarly if $a=c_2$ then $a=\lgf{c_1}<_Qc_1<_Q\rgf{c_1}=\rgf{a}=b$ and one can
take $c=c_1$.
Suppose now that $c_1\not=b$ and $c_2\not=a$. If
$|\StBr _T|(a)\le |\StBr _T|(b)$ then using $a<_Q\rgf{c_2}$ the converse of
(vi) implies that $a\le _Qc_2$. Since $|\StBr _T|(c_2)>|\StBr _T|(a)$ one
obtains $a<_Qc_2$ and hence one can set $c:=c_2$. Similarly one gets
$c:=c_1<_Qb$ if $|\StBr _T|(b)\le |\StBr _T|(a)$.
\end{bew}

Suppose now that $V$ is a \YD module as in Section \ref{s-Nichols}
with $d:=\dim _k V=2$. Let $\{x_1,x_2\}$ denote a canonical basis of $V$.
Then for a full binary tree $T$ the following assignment
defines a unique map $\tau _0:\Nbar (T)\to V^\ot $.
\begin{itemize}
\item
$\tau _0(\lgh ):=x_2$, $\tau _0(\rgh ):=x_1$.
\item
If $a\in N(T)$ then $\tau _0(a):=\tau _0(\rgf{a})\tau _0(\lgf{a})
-\chi (\tau _0(\rgf{a}),\tau _0(\lgf{a}))\tau _0(\lgf{a})\tau _0(\rgf{a})$.
\end{itemize}
Note that one has $\deg (\tau _0(a))=\StBr _T(a)$ for all $a\in \Nbar (T)$.
To shorten notation we will write $\chi (a,b)$ and $g(a)$ instead of
$\chi (\tau _0(a),\tau _0(b))$ and $g(\tau _0(a))$, respectively, for any
$a,b\in \Nbar (T)$.
Let $\tau :\Nbar (T)\to \cB (V)$ denote the composition of $\tau _0$
with the canonical map $V^\ot \to \cB (V)$. 

\begin{defin}
Let $n\in \ndN _0$, $T$ a full binary tree, $V$ a \YD module as in
Section \ref{s-Nichols} with $\dim _kV=2$, and let
$\cB (V)$ denote the corresponding Nichols algebra. We say that
\textit{$\cB (V)$ is of type $T$ in degree $n$} if there exists a
canonical basis $\{x_1,x_2\}$ of $V$ such that for all $a\in \Nbar _2(T)$
with $|\StBr _T|(a)\le n$ the numbers $\chi (a,a)\in k$ are roots of unity
but different from 1 and the sets
\begin{align*}
\bigg\lbrace\prod _{a\in \Nbar _2(T)}\tau (a)^{i_a}\,\Big|\,
0\le i_a<\mathrm{ord}\,\chi (a,a)\quad \forall a\in \Nbar _2(T),\quad
\sum _{a\in \Nbar _2(T)}i_a |\StBr _T|(a)\le m\bigg\rbrace
\end{align*}
form a basis of the $k$ vector spaces $\bigoplus _{i=0}^m\cB (V)_i$,
respectively.
Here the elements of $\Nbar _2(T)$ are ordered with respect to the order
$<_Q$. Further, we say that \textit{$\cB (V)$ is of type $T$} if $\cB (V)$
is of type $T$ in degree $n$ for all $n\in \ndN $.
\end{defin}

Note that $\cB (V)$ is of type $T$ in degree 0 for any full binary tree
$T$. Further, if $\cB (V)$ is of type $T$ for a full binary tree $T$
then $\cB (V)$ is finite dimensional. More exactly, one gets
$\dim _k\cB (V)=\prod _{a\in \Nbar _2(T)}\mathrm{ord}\,\chi (a,a)$.

\subsection{The main result}

Let $R_n$ denote the set of primitive $n^\mathrm{th}$ roots of unity in $k$
where $n\in \ndN $, $n\ge 2$.
The following theorem is the main result of this paper.

\begin{thm}\label{t-class}
\renewcommand{\theenumi }{(T\arabic{enumi})}
Let $k$ be a field of characteristic zero and $G$ an abelian group.
Let $V\in {}_{kG}^{kG}\YDcat $ be a \YD module with $\dim _kV=2$ and completely
reducible $kG$-action. Suppose that there exists a canonical basis of $V$
such that the matrix $(q_{ij})_{i,j=1,2}$ of the braiding of $V$ with respect
to this basis satisfies one of the following conditions.
\begin{enumerate}
\item
$q_{11},q_{22}\in \cup _{n=2}^\infty R_n$, $q_{12}q_{21}=1$.
\item
$(1-q_{11}q_{12}q_{21})(1+q_{11})
=(1-q_{12}q_{21}q_{22})(1+q_{22})=0$,
$q_{12}q_{21}\in \cup _{n=2}^\infty R_n$.
\item
$q_{12}q_{21}=q_{11}^{-2}$, $q_{22}\in \{q_{11}^2,-1\}$,
$q_{11}\in \cup _{n=3}^\infty R_n$, or\\
$q_{11}\in R_3$, $q_{12}q_{21}q_{22}=1$,
$q_{22}\in R_2\cup \cup _{n=4}^\infty R_n$, or\\
$q_{11}\in R_3$, $q_{12}q_{21}=-q_{11}$, $q_{22}=-1$.
\item
$q_0:=q_{11}q_{12}q_{21}\in R_{12}$,
$q_{11}=q_0^4$, $q_{22}=-q_0^2$, or\\
$q_{12}q_{21}\in R_{12}$, $q_{11}=q_{22}=-(q_{12}q_{21})^2$.
\item
$q_{12}q_{21}\in R_{12}$, $q_{11}=-(q_{12}q_{21})^2$,
$q_{22}=-1$, or\\
$q_0:=q_{11}q_{12}q_{21}\in R_{12}$, $q_{11}=q_0^4$, $q_{22}=-1$.
\item
$q_{11}\in R_{18}$, $q_{12}q_{21}=q_{11}^{-2}$, $q_{22}=-q_{11}^3$.
\item
$q_{11}\in R_{12}$, $q_{12}q_{21}=q_{11}^{-3}$, $q_{22}=-1$, or\\
$q_{12}q_{21}\in R_{12}$, $q_{11}=(q_{12}q_{21})^{-3}$, $q_{22}=-1$.
\item
$q_{12}q_{21}=q_{11}^{-3}$, $q_{22}=q_{11}^3$,
$q_{11}\in \cup _{n=4}^\infty R_n$, or\\
$(q_{12}q_{21})^4=-1$, $q_{22}=-1$,
$q_{11}\in \{-q_{12}q_{21},(q_{12}q_{21})^{-2}\}$, or\\
$(q_{12}q_{21})^4=-1$, $q_{11}=(q_{12}q_{21})^2$, $q_{22}=(q_{12}q_{21})^{-1}$.
\item
$q_{12}q_{21}\in R_9$,
$q_{11}=(q_{12}q_{21})^{-3}$, $q_{22}=-1$.
\item
$q_{12}q_{21}\in R_{24}$, $q_{11}=(q_{12}q_{21})^{-6}$,
$q_{22}=(q_{12}q_{21})^{-8}$.
\item
$q_{11}\in R_5\cup R_{20}$, $q_{12}q_{21}=q_{11}^{-3}$,
$q_{22}=-1$.
\item
$q_{11}\in R_{30}$,
$q_{12}q_{21}=q_{11}^{-3}$, $q_{22}=-q_{11}^5$.
\item
$q_{12}q_{21}\in R_{24}$, $q_{11}=(q_{12}q_{21})^6$,
$q_{22}=(q_{12}q_{21})^{-1}$.
\item
$q_{11}\in R_{18}$, $q_{12}q_{21}=q_{11}^{-4}$, $q_{22}=-1$.
\item
$q_{12}q_{21}\in R_{30}$,
$q_{11}=-(q_{12}q_{21})^{-3}$, $q_{22}=(q_{12}q_{21})^{-1}$.
\item
$q_{11}\in R_{10}$, $q_{12}q_{21}=q_{11}^{-4}$, $q_{22}=-1$, or\\
$q_{12}q_{21}\in R_{20}$, $q_{11}=(q_{12}q_{21})^{-4}$,
$q_{22}=-1$.
\item
$q_{12}q_{21}\in R_{24}$,
$q_{11}=-(q_{12}q_{21})^4$, $q_{22}=-1$.
\item
$q_{12}q_{21}\in R_{30}$,
$q_{11}=-(q_{12}q_{21})^5$, $q_{22}=-1$.
\item
$q_{11}\in R_{14}$, $q_{12}q_{21}=q_{11}^{-3}$, $q_{22}=-1$.
\item
$q_{12}q_{21}\in R_{30}$, $q_{11}=(q_{12}q_{21})^{-6}$, $q_{22}=-1$.
\item
$q_{11}\in R_{24}$, $q_{12}q_{21}=q_{11}^{-5}$, $q_{22}=-1$.
\item
$q_{11}\in R_{14}$, $q_{12}q_{21}=q_{11}^{-5}$, $q_{22}=-1$.
\end{enumerate}
Then the Nichols algebra $\cB (V)$ is finite dimensional.
More precisely, the Nichols algebras satisfying the condition in (T$n$)
are of type T$n$ where $n\in \{1,2,\ldots 22\}$ and the full binary tree
T$n$ is given in Appendix \ref{app-types}. 
Further, all relations of $\cB (V)$ are elements of the ideal of
$V^\ot $ generated by the set
\begin{align}\label{eq-B(V)rels}
&\{\tau _0(a)\,|\,a\in N_0(T)\}\cup
\{\tau _0(a)^{\mathrm{ord}\,p_a}\,|\,a\in \Nbar _2(T)\}\cup {}\\
\notag
&\qquad \{\tau _0(b)\tau _0(\lgf{c})-\chi (b,\lgf{c})\tau _0(\lgf{c})\tau _0(b)
-\mu (b)/\qfact{\rgfl (b)+1}{p_c}\tau _0(c)^{\rgfl (b)+1}\,|\\
\notag
&\phantom{\qquad \{\tau _0(b)\tau _0(\lgf{c})-\chi (b,\lgf{c})
\tau _0(\lgf{c})\tau _0(b)-\mu (b)/}
b\in N_2(T),c:=\lgf{b}\in N_2(T)\}.
\end{align}
\end{thm}

The proof of Theorem \ref{t-class} will be given at the end of this paper.

\begin{bem}
There are already a lot of rank 2 Nichols algebras which are known to
be finite dimensional. Those with $q_{12}q_{21}=q_{11}^{-n}$,
$q_{22}=q_{11}^n$,
$n\in \{0,1,2,3\}$, are called \textit{of finite Cartan type} in
\cite{inp-AndrSchn02}. They appear in (T1), (T2), (T3), and (T8).
Other examples
(cf.~Section 3.3 in \cite{inp-AndrSchn02}) cover essentially all of (T2) and
(T3). Further, there exist recent computations on Nichols algebras by
M.~Gra\~na and Ch. Heaton \cite{privat-Grana}
which give a PBW basis for the examples in (T5) and (T9).
\end{bem}

\section{Finiteness of the Nichols algebras}
\label{s-finiteness}

\subsection{Lyndon words and full binary trees}
\label{ss-Lyndonfull}

In \cite{a-Khar99} Kharchenko proves that any finite dimensional Nichols
algebra
of diagonal type has a Poincar\'e--Birkhoff--Witt basis. He also
gives very useful information on the relations of the algebra. It is worth to
mention that (as was shown by
Ufer in \cite{a-Ufer03}) such results hold in a more general context, namely
for
Nichols algebras generated by a braided vector space with triangular braiding.
In order to prove finite dimensionality of the Nichols algebras in Theorem 6
results from \cite{b-Loth83} and \cite{a-Khar99} are recalled and adapted to
our
conventions. We consider only the rank 2 case and replace the symbol $>$ for
the lexicographic order in \cite{a-Khar99} by $<$.

Set $X:=\{\alp ,\bet \}$ and consider the total order $<$ on $X$ given by
$\alp <\bet $. Let $\mbbX $ and $\mbbX ^+$ denote the set of words and
nonempty words, respectively, in the letters $\alp ,\bet $. Then $<$ induces
the lexicographic order on $\mbbX $: $u,v\in \mbbX $ satisfy $u<v$ if
and only if either $v=uw$ for some $w\in \mbbX ^+$ or there exist
$u_1,u_2,u_3\in \mbbX $ such that $u=u_1\alp u_2$ and $v=u_1\bet u_3$.
For $u,v\in \mbbX $ we write $u\le v$ if $u=v$ or $u<v$.
The length of a word $u$, i.\,e.~the number of its letters, is denoted by
$|u|$. A word $u\in \mbbX ^+$ is called \textit{Lyndon}
if for any decomposition $u=vw$ with $v,w\in \mbbX ^+$ the relation
$vw<wv$ holds.

\begin{satz}\label{s-Lyndon}
(i) \cite[Prop.~5.1.2]{b-Loth83} A word $u\in \mbbX ^+$ is Lyndon
if and only if $u=vw$ with $v,w\in \mbbX ^+$ implies $u<w$.\\
(ii) \cite[Prop.~5.1.3]{b-Loth83} A word $u\in \mbbX ^+$ is Lyndon
if and only if either $u\in X$ or there exist Lyndon words $v,w\in \mbbX $
such that $v<w$ and $u=vw$.
\end{satz}

Any word $u\in \mbbX $ has a unique decomposition into the product of a
nonincreasing sequence of Lyndon words \cite[Thm.~5.1.5]{b-Loth83}.
Further, any Lyndon word $u\notin X$ has
a decomposition into the product of two Lyndon words $u=vw$ (which then
satisfy $v<w$) such that $|v|$ is minimal.
This is called the \textit{Shirshow decomposition} of $u$.

\begin{satz}\cite[Prop.~5.1.4]{b-Loth83}\label{s-Shirshow}
Suppose that $u,v,w$ are Lyndon words and $u=vw$. Then $u=vw$ is the Shirshow
decomposition of $u$ if and only if $v\in X$ or for the Shirshow
decomposition $v=v_1v_2$ the relation $w\le v_2$ holds.
\end{satz}

Let $k$ be a field of characteristic zero, $G$ an abelian group, and
$V\in {}^{kG}_{kG}\YDcat $ a two dimensional \YD module of diagonal type.
Let $\cI $ be an $\ndN _0^2$-graded ideal of $V^\ot $, $\cA $ the
$\ndN _0^2$-graded algebra $V^\ot /\cI $, and $\cA _n$ the subspace of
$\cA $ of homogeneous elements of total degree $n$
(where $n\in \ndN _0$). Let $\cA ^+$ denote the unique maximal ideal
of $\cA $. The following results will be needed for
$\cA =V^\ot $ and for $\cA =\cB (V)$.

After fixing a canonical basis $\{x_2,x_1\}$ of $V$ one can associate to any
Lyndon word $u\in \mbbX $ an element
$\suplet{u}\in \cA $ as follows. Set $\suplet{\alp }:=x_2$,
$\suplet{\bet }:=x_1$, and $\suplet{u}:=\suplet{w}\suplet{v}-\chi (\suplet{w},
\suplet{v})\suplet{v}\suplet{w}$ if $u=vw$ is the Shirshow decomposition of
$u$.
Note that this definition differs from that in \cite{a-Khar99} by a constant
for each Lyndon word $u$. However this is not relevant for the following
assertions.

\begin{lemma}\cite[Lemma~3]{a-Khar99}\label{l-uhv}
If $u,v$ are Lyndon words with $u<v$ then $u^h<v$ for any $h>0$.
\end{lemma}

\begin{lemma}\cite[Lemma~5]{a-Khar99}\label{l-[u]poly}
Let $u$ be a Lyndon word with $|u|=m$. Then $\suplet{u}\in \cA $
is a linear combination of monomials $\suplet{\alpha _1}\suplet{\alpha _2}
\cdots \suplet{\alpha _m}\in \cA $, $\alpha _i\in X$, such that
$u\le \alpha _1\alpha _2\cdots \alpha _m$.
\end{lemma}

\begin{lemma}\cite[Lemma~6]{a-Khar99}\label{l-freerels}
If $u<v\in \mbbX $ are Lyndon words then $\suplet{v}\suplet{u}-
\chi (\suplet{v},\suplet{u})\suplet{u}\suplet{v}$
is a $k$ linear combination of products
$\suplet{u_1}\suplet{u_2}\cdots \suplet{u_i}$ for certain $i\in \ndN $ and
Lyndon words $u_j$ with $u<u_j<v$ such that $\deg (\suplet{v}\suplet{u})=
\deg (\suplet{u_1}\suplet{u_2}\cdots \suplet{u_i})$ and
$uv\le u_1u_2\cdots u_i$.
\end{lemma}

For a Lyndon word $u\in \mbbX $ let $\hght{u}\in \ndN $ denote the smallest
number such that $\suplet{u}^{\hght{u}}$ can be written as a linear
combination of products $\suplet{u_1}\suplet{u_2}\cdots \suplet{u_i}$,
$i\in \ndN $, where $u_j$ are Lyndon words with $u<u_j$. By Lemma \ref{l-uhv}
the relation $u<u_1$ implies that $u^{\hght{u}}<u_1$. Now since
$u^{\hght{u}}$ is not the beginning of $u_1$ one obtains that
$u^{\hght{u}}<u_1u_2\ldots u_i$ has to hold as well.

Define $B:=\{u\in \mbbX \,|\,u\mbox{ is Lyndon},\hght{u}>1\}$.
For each $u\in B$ let $S(u)_<$ and $S(u)$
denote the subalgebras of $\cA $ generated by the sets
$\{\suplet{v}\,|\,v\in B,u<v\}$
and $\{\suplet{v}\,|\,v\in B,u\le v\}$, respectively.
Define $S(u)^+:=S(u)\cap \cA ^+$ and
$S(u)_<^+:=S(u)_<\cap \cA ^+$.

\begin{thm}\cite[Theorem~2]{a-Khar99}\label{t-BVbasis}
The set $\{\suplet{u_1}^{n_1}\suplet{u_2}^{n_2}\cdots \suplet{u_i}^{n_i}\,|\,
i\in \ndN _0,u_j\in B,u_i<\cdots <u_2<u_1,n_j<\hght{u_j}\,\forall j\}$
is a basis of the $k$ vector space $\cA $.
\end{thm}

\begin{folg}\label{f-S(u)basis}
For any $n\in \ndN $, $u\in B$ the sets
\begin{align*}
\{\suplet{u_1}^{n_1}\suplet{u_2}^{n_2}
\cdots \suplet{u_i}^{n_i}\,|\,i\in \ndN _0,u_j\in B,
u<u_1<u_2<\cdots <u_i,n_j<\hght{u_j}\,\forall j\}\\
\intertext{and}
\{\suplet{u_1}^{n_1}\suplet{u_2}^{n_2}
\cdots \suplet{u_i}^{n_i}\,|\,i\in \ndN _0,u_j\in B,
u\le u_1<u_2<\cdots <u_i,n_j<\hght{u_j}\,\forall j\}
\end{align*}
form a basis of the $k$ vector space $S(u)_<$ and $S(u)$, respectively.
\end{folg}

\begin{bew}
Since $S(u)_<=\bigcup _{v\in B,u<v}S(v)$ it suffices to prove the assertion
for $S(u)$.
As $\cA =\bigoplus _{n=0}^\infty \cA _n$ and
$S(u)=\bigoplus _{n=0}^\infty S(u)\cap \cA _n$
the proof can be performed by induction on $n$.
Note that for given $n\in \ndN _0$ the set $\cA _n$
is finite dimensional and $\{v\in B\,|\,\deg (\suplet{v})\le n\}$
is a finite set. Suppose that $\{u_1,u_2,\ldots ,u_i\}=\{v\in B\,
|\,u<v,\deg (\suplet{v})\le n\}$ and $u_1<u_2<\ldots <u_i$.
By Lemma \ref{l-freerels} the elements
$\suplet{u}^m\suplet{u_{j_1}}\suplet{u_{j_2}}\cdots \suplet{u_{j_r}}$,
$m<\hght{u}$, $m\deg (\suplet{u})+\sum _{s=1}^r\deg (\suplet{u_{j_s}})=n$,
span $S(u)\cap \cA _n$. By induction hypothesis the elements
$\suplet{u_1}^{n_1}\suplet{u_2}^{n_2}\cdots \suplet{u_i}^{n_i}$,
$n_j<\hght{u_j}$ for all $j$, where $\sum _{j=1}^in_j\deg (\suplet{u_j})=
n-m\deg (\suplet{u})$, span $S(u_1)\cap \cA _{n-m\deg (\suplet{u})}$ for
all $m>0$. Therefore the elements
$\suplet{u}^m\suplet{u_1}^{n_1}\suplet{u_2}^{n_2}\cdots \suplet{u_i}^{n_i}$,
$n_j<\hght{u_j}$ for all $j$, $m<\hght{u}$, where
$m\deg(\suplet{u})+\sum _{j=1}^in_j\deg (\suplet{u_j})=n$, span
$S(u)\cap \cA _n$. Then Theorem
\ref{t-BVbasis} and a simple dimension argument imply the assertion.
\end{bew}

Later we will need the fact that for $w\in B$ one has
\begin{align}\label{eq-S(u)zerl}
S(w)=(S(w)S(w)_<^+ +S(w)\suplet{w}^2)\oplus k\suplet{w}\oplus k1
\end{align}
as graded vector spaces which is one of the consequences of
Corollary \ref{f-S(u)basis} and Lemma \ref{l-freerels}.

\begin{lemma}\cite[Corollary 2]{a-Khar99}\label{l-hght1}
For a Lyndon word $u$ of length $l_u$ one has $\hght{u}=1$ if and only if
$\suplet{u}$ can be written as a linear combination of monomials
$\suplet{\alpha _1}\suplet{\alpha _2}\cdots \suplet{\alpha _{l_u}}$,
$\alpha _i\in X$ for all $i$,
such that $u<\alpha _1\alpha _2\cdots \alpha _{l_u}$.
\end{lemma}

\begin{lemma}\label{l-tree2Lyndon}
Let $T$ be a full binary tree.\\
(i) There exists a unique map $\gamma :\Nbar (T)
\to \mbbX $ such that $\gamma (\lgh )=\alp $, $\gamma (\rgh )=\bet $, and
$\gamma (a)=\gamma (\lgf{a})\gamma (\rgf{a})$ for all $a\in N(T)$.\\
(ii) For $a\in \Nbar (T)$ the equation $|\gamma (a)|=|\StBr _T|(a)$
holds.\\
(iii) Any word $\gamma (a)$, $a\in \Nbar (T)$, is Lyndon and
$\gamma (\lgf{a})\,\gamma (\rgf{a})$ is the Shirshow decomposition of
$\gamma (a)$ for $a\in N(T)$.\\
(iv) For any $a,b\in \Nbar (T)$ the relation
$\gamma (a)<\gamma (b)$ is equivalent to $a<_Qb$.
\end{lemma}

\begin{bew}
(i) Existence and uniqueness of $\gamma $ follow from the facts that
$\lgf{(\cdot )}$ and $\rgf{(\cdot )}$ are well defined maps from $N(T)$ to
$\Nbar (T)$ and $|\StBr _T|(\lgf{a})<|\StBr _T|(a)$, $|\StBr _T|(\rgf{a})<
|\StBr _T|(a)$ for all $a\in N(T)$.

(ii) This follows immediately from the definition of $\gamma $ and
$|\StBr _T|$.

(iii), (iv) We use induction on $|\StBr _T|(a)$ and
$\max \{|\StBr _T|(a),|\StBr _T|(b)\}$,
respectively. If $a=\lgh $ or $a=\rgh $ then $\gamma (a)$ is Lyndon. Further,
if $|\StBr _T|(a)=|\StBr _T|(b)=1$ then $a,b\in \{\lgh ,\rgh \}$ and hence
$a<_Qb$ is equivalent to $a=\lgh ,b=\rgh $ which holds if and only if
$\gamma (a)=\alp <\bet  =\gamma (b)$.

Assume now that (iii) and (iv) hold
whenever $a,b\in \Nbar (T)$, $|\StBr _T|(a)\le n$, and $|\StBr _T|(b)\le n$
for some $n\in \ndN $. If $a\in N(T)$ then by induction hypothesis
$\gamma (\lgf{a})$ and $\gamma (\rgf{a})$ are Lyndon words.
Since $\lgf{a}<_Q\rgf{a}$ we also have $\gamma (\lgf{a})<\gamma (\rgf{a})$.
Thus $\gamma (\lgf{a})\gamma (\rgf{a})$ is Lyndon by Proposition
\ref{s-Lyndon}(ii). This proves the induction step of the first part of (iii).

Now we prove (iv) in the case $|\StBr _T|(a)= n+1$, $|\StBr _T|(b)\le n$. The
proof for $|\StBr _T|(b)=n+1$ is completely analogous.
Let $(a_1,a_2,\ldots ,a_m)$ denote the set of nodes of $T$ with
$|\StBr _T|(a_i)\le n$ in increasing order with respect to $<_Q$.
By Lemma \ref{l-fbtreeord}(ix) the node $a\in N(T)$
is the unique $c\in \Nbar (T)$ such that $|\StBr _T|(c)\le n+1$ and
$\lgf{a}<_Qc<_Q\rgf{a}$. Thus there exists $i\in \ndN $ such that
$\lgf{a}=a_i$ and $\rgf{a}=a_{i+1}$. On the other hand,
the induction hypothesis gives that $\gamma (a_j)<\gamma (a_l)$ if and only if
$j<l$. Now note that $\gamma (a_i)<\gamma (a)=\gamma (a_i)\gamma (a_{i+1})<
\gamma (a_{i+1})$ as $\gamma (a)$ is Lyndon.

It remains to show that $\gamma (\lgf{a})\gamma (\rgf{a})$ is the Shirshow
decomposition of $\gamma (a)$ where $|\StBr _T|(a)=n+1$. If
$\lgf{a}=\lgh $ then $|\gamma (\lgf{a})|=1$ and we are
done. Otherwise $\lgf{a}\in N(T)$ and Lemma \ref{l-fbtreeord}(vi),(vii)
gives $\rgf{a}\le _Q\rgf{\lgf{a}}$. Therefore $\gamma (\rgf{a})\le
\gamma (\rgf{\lgf{a}})$ by the induction hypothesis for (iv).
Further, the induction hypothesis of (iii) tells us that
$\gamma (\lgf{\lgf{a}})\gamma (\rgf{\lgf{a}})$ is the Shirshow
decomposition of $\gamma (\lgf{a})$.
Thus Proposition \ref{s-Shirshow} for $u=\gamma (a)$ together with the last two
relations give the claim.
\end{bew}

Our aim in this section is to give a computable criterion which ensures
that the Nichols algebras in Theorem \ref{t-class} are finite dimensional
and of the given type. To do so we have to introduce additional notation
which will be needed only for $\cA =\cB (V)$.

For a Lyndon word $u$ and $n\in \ndN $ let $\cF (u)_n$ denote the
$kG$-module
\begin{align*}
\cF (u)_n:=\bigg(S(u)\cap \big(\cB (V)S(u)_<^+
+\cB (V)\suplet{u}^2\big) \cap \bigoplus _{m=1}^n\cB (V)_m\bigg)\# kG
\end{align*}
and set $\cF (u)_0=\{0\}$. By Corollary \ref{f-S(u)basis} one obtains that
\begin{align*}
\cF (u)_n=\bigg(\big(S(u)S(u)_<^++S(u)\suplet{u}^2\big) \cap
\bigoplus _{m=1}^n\cB (V)_m\bigg)\# kG.
\end{align*}

Let $n\in \ndN $. Suppose that there exists a full binary tree $T$
such that for any Lyndon word $u$ with $|u|\le n$ the relation
$\hght{u}>1$ is equivalent to $u=\gamma (a)$ for some $a\in \Nbar _2(T)$.
The definition of $\tau $ and $\suplet{\cdot }$ and Lemma \ref{l-tree2Lyndon}(iii)
imply that for any $a\in \Nbar _2(T)$ one has $\tau (a)=\suplet{\gamma (a)}$.
Then by Corollary \ref{f-S(u)basis} with $u=\alp $ and by Lemma
\ref{l-tree2Lyndon}(iv) the set
\begin{align*}
\bigg\{\prod _{a\in \Nbar _2(T)}\tau (a)^{i_a}\,\bigg|\,i_a<\hght{\gamma (a)},
\sum _{a\in \Nbar _2(T)}i_a\mathrm{totdeg}\,(\tau (a))\le n\bigg\}
\end{align*}
where the product is taken with respect to the order $<_Q$ of $\Nbar _2(T)$
forms a basis of $\bigoplus _{m=0}^n\cB (V)_m$.

For $a\in \Nbar (T)$ define $p_a:=\chi (a,a)^{-1}$ and
\begin{align}\label{eq-lambda}
\lambda (a):=\begin{cases}
0 & \text{if $a\notin N(T)$,}\\
\chi (\lgh ,\rgh )^{-1}-\chi (\rgh ,\lgh ) & \text{if $a=\roo $,}\\
\chi (\lgf{a},\rgf{a})^{-1}-\chi (\rgf{a},\lgf{a})+\lambda (b)
& \text{otherwise,}
\end{cases}
\end{align}
where $b$ is the parent of $a$. Further, for any $b\in N(T)$ with
$\lgf{b}\in N(T)$ set
\begin{align*}
\mu (b):=&\begin{cases}
\lambda (b) & \text{if $b=\rch{\lgf{b}}$,}\\
\lambda (b)\mu (\rgf{b}) & \text{otherwise.}
\end{cases}
\end{align*}
Finally, for any $b\in N(T)$ with $c:=\lgf{b}\in N(T)$ and
$\rgfl (b)\le 2$ set
\begin{align*}
\nu (b):=\begin{cases}
\chi (\lgf{c},b)^{-1}-\chi (b,\lgf{c})
+\lambda (b)\lambda (c)\Big(\qnum{2}{p_f}^{-1}-\qnum{2}{p_c}^{-1}\Big) &
\text{if $\rgfl (b)=1$,}\\
\chi (\lgf{c},\rch{c})^{-1}+\lambda (c)\lambda (\rch{c})\qnum{2}{p_c}^{-1}
\Big(\qnum{2}{p_f}^{-1}-\qnum{3}{p_c}^{-1}\Big) & \text{if $\rgfl (a)=2$,}
\end{cases}
\end{align*}
where $f=\rgf{c}$,
whenever all denominators are nonzero.

\subsection{The finiteness results}
\label{ss-finiteness}

{}From now on let $T$ be a full binary tree such that
\begin{align}\label{eq-inequality}
\min \{\rchl (\lch{b}),\lchl (\rch{b})\}\le 3\quad
\text{for all $b\in N_2(T)$},
\end{align}
i.\,e.~either $\rch{\rch{\rch{\lch{b}}}}\notin N(T)$ or
$\lch{\lch{\lch{\rch{b}}}}\notin N(T)$.
Note that all binary trees in Appendix \ref{app-types} satisfy
this condition.

\begin{defin}
We call a triple $(T,V,n)$ where $n\in \ndN _0$
and $V\in {}^{kG}_{kG}\YDcat $ is a two-dimensional \YD module
of diagonal type \textit{admissible} if the following conditions
are satisfied.
\begin{itemize}
\item
For $a\in N(T)$ with $|\StBr _T|(a)\le n$ relation $a\in N_2(T)$
is equivalent to $\lambda (a)\not=0$,
\item
the numbers $p_a$ for $a\in \Nbar _2(T)$, $|\StBr _T|(a)\le n$,
are roots of unity but different from 1,
\item
if $a\in N_2(T)$, $\lch{a}\in N_2(T)$, and $|\StBr _T|(\lch{a})\le n$
then $p_a\not=-1$, $p_{\rgf{a}}\not=-1$,
\item
for all $b\in N_2(T)$ with $c:=\lgf{b}\in N_2(T)$ and
$|\StBr _T|(b)+|\StBr _T|(\lgf{\lgf{b}})\le n$ one has
$\qfact{\rgfl (b)+1}{p_c}\not=0$ and either
$\rchl (\lch{c})\le \rgfl (b)$ or $\rgfl (b)\le 2$, $\nu (b)=0$.
\end{itemize}
\end{defin}

Note that if $b\in N_2(T)$ and $c:=\lgf{b}\in N_2(T)$ then
$\rgfl (b)<\lchl (\rch{c})$. Thus by
(\ref{eq-inequality}) for $c$ one has either $\rgfl (b)\le 2$ or
$\rchl (\lch{c})\le \rgfl (b)$.

\begin{satz}\label{s-findegn}
Let $T$ denote a full binary tree and $V\in {}^{kG}_{kG}\YDcat $ a
two-dimensional \YD module of diagonal type. Let $n\in \ndN _0$. If
$(T,V,n)$ is admissible then $\cB (V)$ is of type $T$ in degree $n$
and all homogeneous relations of $\cB (V)$ of degree at most $n$ are
elements of the ideal of $V^\ot $ generated by the set
\begin{align*}
&\{\tau _0(a)\,|\,a\in N_0(T),|\StBr _T|(a)\le n\}\cup {}\\ 
&\qquad \{\tau _0(a)^{\mathrm{ord}\,p_a}\,|\,a\in \Nbar _2(T),|\StBr _T|(a)
\cdot \mathrm{ord}\,p_a\le n\}\cup {}\\
&\qquad \{\tau _0(b)\tau _0(\lgf{c})-\chi (b,\lgf{c})\tau _0(\lgf{c})\tau _0(b)
-\mu (b)/\qfact{\rgfl (b)+1}{p_c}\tau _0(c)^{\rgfl (b)+1}
\,|\\
&\qquad \qquad \qquad \qquad b\in N_2(T),c:=\lgf{b}\in N_2(T),
|\StBr _T|(b)+|\StBr _T|(\lgf{c})\le n\}.
\end{align*}
\end{satz}

The following corollary is an immediate consequence of Proposition
\ref{s-findegn}. It will be our main tool to prove Theorem \ref{t-class}.

\begin{folg}\label{f-findeg}
Let $T$ denote a full binary tree and $V\in {}^{kG}_{kG}\YDcat $ a
two-dimensional \YD module of diagonal type. If
$(T,V,n)$ is admissible for all $n\in \ndN $ then $\cB (V)$ is of type $T$
and all relations of $\cB (V)$ are
elements of the ideal of $V^\ot $ generated by the set
\begin{align*}
&\{\tau _0(a)\,|\,a\in N_0(T)\}\cup
\{\tau _0(a)^{\mathrm{ord}\,p_a}\,|\,a\in \Nbar _2(T)\}\cup {}\\
&\qquad \{\tau _0(b)\tau _0(\lgf{c})-\chi (b,\lgf{c})\tau _0(\lgf{c})\tau _0(b)
-\mu (b)/\qfact{\rgfl (b)+1}{p_c}\tau _0(c)^{\rgfl (b)+1}\,|\\
&\phantom{\qquad \{\tau _0(b)\tau _0(\lgf{c})-\chi (b,\lgf{c})
\tau _0(\lgf{c})\tau _0(b)-\mu (b)/}
b\in N_2(T),c:=\lgf{b}\in N_2(T)\}.
\end{align*}
\end{folg}

\begin{bew}[ of Proposition \ref{s-findegn}]
We proceed by induction over $n$. As noted previously the assertion is
true for $n=0$. Assume that Proposition \ref{s-findegn} is valid for
$(T,V,n-1)$ and that $(T,V,n)$ is admissible.
By Corollary \ref{f-S(u)basis} for $\cA =V^\ot $ and with $u=\alp $ and
by Lemma \ref{l-tree2Lyndon} it suffices to prove that the following
assertions hold.
\begin{itemize}
\item[(a)]
If $a\in \Nbar _2(T)$ and $|\StBr _T|(a)\le n$ then
$\hght{\gamma (a)}\ge \min \{\mathrm{ord}\,p_a,(n+1)/|\StBr _T|(a)\}$. If
$n\ge |\StBr _T|(a)\cdot \mathrm{ord}\,p_a$ then
$\suplet{\gamma (a)}^{\mathrm{ord}\,p_a}=0$.
\item[(b)]
If $a\in N_0(T)$ and $|\StBr _T|(a)\le n$ then $\suplet{\gamma (a)}=0$.
\item[(c)]
If $u$ is Lyndon, $2\le |u|\le n$, and $u=\gamma (a)\gamma (b)$ is the
Shirshow decomposition of $u$ with
$a,b\in \Nbar _2(T)$, $a=\lgf{c}$, where $c=\lgf{b}$, then
$\suplet{u}=\mu (b)/\qfact{\rgfl (b)+1}{p_c}\tau (c)^{\rgfl (b)+1}$ and
$u<\gamma (c)$.
\item[(d)]
If $u$ is Lyndon, $|u|\le n$, $u\notin \{\gamma (a)\,|\,a\in \Nbar (T)\}$,
and $u$ is not as in (c) then $\hght{u}=1$
and the relation corresponding to $u$ follows from those given in
Proposition \ref{s-findegn} for $(T,V,n-1)$. 
\end{itemize}
In order to prove assertions (a)--(d) we additionally
use the following induction hypotheses which will be proven after the proof
of (a)--(d).
\begin{itemize}
\item[(e)]
If $a\in N(T)$, $m:=|\StBr _T|(a)\le n$, and $u=\gamma (a)$ has Shirshow
decomposition $u=vw$ then
\begin{align}\label{eq-copr[u]}
\begin{aligned}
\copr (\iota (\suplet{u}))&-\iota (\suplet{u})\ot 1-g(\suplet{u})^{-1}\ot
\iota (\suplet{u})\\
&-\lambda (a)g(\suplet{v})^{-1}\iota (\suplet{w})\ot \iota (\suplet{v})\in
\iota (\cF (w)_{m-1})\ot \cB (V^*).
\end{aligned}
\end{align}
\item[(f)]
If $a,b\in \Nbar (T)$ such that $|\StBr _T|(a)\le n$, $|\StBr _T|(b)\le n$,
and $a<_Qb$ then $\paar{\iota (\tau (b))}{\tau (a)}=0$.
\item[(g)]
If $a\in N(T)$ and $|\StBr _T|(a)\le n$ then the following equations hold.
\begin{align*}
\paar{\iota (\tau (\lgf{a}))}{\tau (a)}
=&\lambda (a)\paar{\iota (\tau (\lgf{a}))}{\tau (\lgf{a})}\tau (\rgf{a}),\\
\paar{\iota (\tau (a))}{\tau (a)}
=&\lambda (a)\paar{\iota (\tau (\lgf{a}))}{\tau (\lgf{a})}
\paar{\iota (\tau (\rgf{a}))}{\tau (\rgf{a})}.
\end{align*}
\item[(h)]
Suppose that $b\in N_2(T)$, $c:=\lgf{b}\in N(T)$, and $|\StBr _T|(b)+
|\StBr _T|(\lgf{c})\le n$. Set $a:=\lgf{c}$.
If $\rgfl (b)=1$ then $b=\rch{c}$ (see Figure 1) and one has
\begin{align*}
\paarb{\iota (\tau (a))}{\tau (b)\tau (a)-\chi (b,a)
\tau (a)\tau (b)-\frac{\lambda (b)}{\qnum{2}{p_c}}\tau (c)^2}\qquad &\\
=\spaar{a}{a}\nu (b)\tau (b)&.
\end{align*}
If $\rgfl (b)=2$ then $b=\lch{\rch{c}}$ and with $d:=\lch{c}$
(see Figure 2) one has
\begin{align*}
\paarb{\iota (\tau (d))}{\tau (b)\tau (a)-\chi (b,a)
\tau (a)\tau (b)-\frac{\mu (b)}{\qfact{3}{p_c}}\tau (c)^3}\qquad &\\
=\lambda (b)\spaar{d}{d}\nu (b)\tau (\rch{c})&.
\end{align*}
\end{itemize}
\begin{center}
\setlength{\unitlength}{4pt}
\hfill
\begin{picture}(12,8)
\put(0,4){\line(1,1){4}}
\put(4,8){\line(1,-1){8}}
\put(4,0){\line(1,1){4}}
\put(4,8){\circle*{1}}
\put(0,4){\circle*{1}}
\put(8,4){\circle*{1}}
\put(4,0){\circle*{1}}
\put(12,0){\circle*{1}}
\put(0,8){$a$}
\put(3,4){$c$}
\put(7,0){$b$}
\end{picture}\hfill
\begin{picture}(16,12)
\put(0,4){\line(1,1){8}}
\put(4,8){\line(1,-2){2}}
\put(8,12){\line(1,-1){8}}
\put(12,8){\line(-1,-2){4}}
\put(10,4){\line(1,-2){2}}
\put(8,12){\circle*{1}}
\put(4,8){\circle*{1}}
\put(12,8){\circle*{1}}
\put(0,4){\circle*{1}}
\put(6,4){\circle*{1}}
\put(10,4){\circle*{1}}
\put(16,4){\circle*{1}}
\put(8,0){\circle*{1}}
\put(12,0){\circle*{1}}
\put(3,4){$d$}
\put(1,10){$a$}
\put(7,9){$c$}
\put(13,10){$f$}
\put(11,4){$\rch{c}$}
\put(9,0){$b$}
\end{picture}
\hfill \makebox[0pt]{}\\
\hfill
Figure 1\footnotemark[1]
\footnotetext[1]{The labels belong to the nodes above them.}
\hfill
Figure 2\footnotemark[1]
\hfill \makebox[0pt]{}
\end{center}

Note that for $n=0$ all assertions (a)--(h) are trivially fulfilled and
hence we may start with the induction step.

To (d).
Suppose that $u$ is a Lyndon word with $|u|=n\ge 2$ and Shirshow
decomposition $u=vw$ and $u$ is as in (d).
If $\hght{v}=1$ or $\hght{w}=1$ then $\hght{u}=1$ by Lemma
\ref{l-[u]poly} and Lemma \ref{l-hght1}. More precisely, in order to ensure
$\hght{u}=1$ one does not need a new relation. 
Otherwise by induction hypothesis there
exist $a,b\in \Nbar _2(T)$ such that $\gamma (a)=v$, $\gamma (b)=w$.
Note that since $\StBr _T(a)=\deg (\suplet{v})$ and
$\StBr _T(b)=\deg (\suplet{w})$ Lemma
\ref{l-fbtreeord}(iii) implies that $a$ and $b$ are uniquely determined by
$\gamma (a)$ and $\gamma (b)$, respectively.

Since $u$ is Lyndon we have also $v<w$ and hence $a<_Qb$ by Lemma
\ref{l-tree2Lyndon}(iv). Further,
Proposition \ref{s-Shirshow} and Lemma \ref{l-tree2Lyndon}(iv) imply that
either $a=\lgh $ or $a\in N_2(T)$ and $b\le _Q\rgf{a}$. Assume for a moment
that $a\in N_2(T)$ and $b=\rgf{a}$ and set $c:=\rch{a}$.
Then $a=\lgf{c}$, $b=\rgf{c}$, and $|\StBr _T|(c)=n$ and hence $u=\gamma (c)$.
This case is covered by (a) and (b).
Similarly, if $a=\lgh $ and $b=\rgh $ then $n=2$ and $u=\gamma (\roo )$
which again belongs to (a) or (b). Therefore one has either $a=\lgh $
and $b\in N_2(T)$ or $a,b\in N_2(T)$ and $a<_Qb<_Q\rgf{a}$.

Now we prove that $a<_Q\lgf{b}$. First note that since $a<_Qb$ Lemma
\ref{l-fbtreeord}(vi) implies the relation $|\StBr _T|(a)<|\StBr _T|(b)$
in the case $a\in N_2(T)$. If $a=\lgh $ then this relation is trivial.
Therefore one gets $a\le _Q\lgf{b}$ by Lemma \ref{l-fbtreeord}(v).
Assume for a moment that $a=\lgf{b}$ and set $c:=\lch{b}$.
Then $a=\lgf{c}$ and $b=\rgf{c}$ and hence $u=\gamma (c)$ which is
covered by (a) and (b).
Thus we are arrived at the situation that $b\in N_2(T)$, $a<_Q\lgf{b}$,
and either $a=\lgh $ or $a\in N_2(T)$ and $b<_Q\rgf{a}$.

If $a=\lgh $ then $a<_Q\lgf{b}$ implies that
$|\StBr _T|(a)<|\StBr _T|(\lgf{b})$. If $a\in N_2(T)$ then the same relation
follows from $a<_Q\lgf{b}$ and Lemma \ref{l-fbtreeord}(vi) as
$\lgf{b}<_Qb<_Q\rgf{a}$. Therefore Lemma \ref{l-fbtreeord}(v) gives
$a\le _Q\lgf{\lgf{b}}$. The case $a=\lgf{\lgf{b}}$ can be omitted as it
is exactly the situation in (c).
Otherwise set $u_1=\gamma (a)\gamma (\lgf{b})$. Note that $u_1$ is
Lyndon by Proposition \ref{s-Lyndon}(ii). We show that
$\hght{u_1}=1$ holds which proves (d) by Lemma \ref{l-[u]poly}
and Lemma \ref{l-hght1}.

Suppose that $\hght{u_1}>1$. As the length of $u_1$ is less than $n$ but at
least 2 by induction hypothesis (a) there exists $c\in N_2(T)$ such that
$u_1=\gamma (c)$.
Since $\gamma (a)<u_1<\gamma (\lgf{b})$ one obtains $a<_Qc<_Q\lgf{b}$.
By Lemma \ref{l-fbtreeord}(iv),(v) this implies that $a\le _Q\lgf{c}$. As
$\gamma (\lgf{c})\gamma (\rgf{c})$ is the Shirshow decomposition of
$\gamma (c)$ by Lemma \ref{l-tree2Lyndon}(iii) the equation
$\gamma (a)\gamma (\lgf{b})=\gamma (\lgf{c})\gamma (\rgf{c})$ implies that
$|\gamma(\lgf{b})|\le |\gamma (\rgf{c})|$ and hence $\gamma (\lgf{c})\le
\gamma (a)$. Thus $a\le _Q\lgf{c}$, Lemma \ref{l-tree2Lyndon}(iv), and
Lemma \ref{l-fbtreeord}(iii) imply that $\lgf{c}=a$ and $\rgf{c}=\lgf{b}$.
Now we have $\lgf{c}=a<_Q\lgf{\lgf{b}}<_Q\lgf{b}=\rgf{c}$, $|\StBr _T|(
\lgf{\lgf{b}})<|\StBr _T|(\lgf{b})$, and $|\StBr _T|(c)>|\StBr _T|(\lgf{b})$.
Then Lemma \ref{l-fbtreeord}(ix) applied to the pair $(c,\lgf{\lgf{b}})$
gives a contradiction.

To (c). The relation $\gamma (a)\gamma (b)<\gamma (c)$ follows from
$\gamma (c)=\gamma (a)\gamma (\rgf{c})$ and $b<_Q\rgf{\lgf{b}}=\rgf{c}$
(see Lemma \ref{l-fbtreeord}(vii)).

Set $U:=\tau (b)\tau (a)-\chi (b,a)\tau (a)\tau (b)-
\mu (b)/\qfact{\rgfl (b)+1}{p_c}\tau (c)^{\rgfl (b)+1}$.
First note that using (e), (f), and (g) one obtains
\begin{align*}
&\paar{\iota (\tau (c))}{\tau (b)\tau (a)-\chi (b,a)\tau (a)\tau (b)-
\mu (b)/\qfact{\rgfl (b)+1}{p_c}\tau (c)^{\rgfl (b)+1}}=\\
&\quad =\paar{\iota (\tau (c))}{\tau (b)}\tau (a)-\chi (\rgf{b},a)\tau (a)
\paar{\iota (\tau (c))}{\tau (b)}\\
&\qquad -(\mu (b)/\qfact{\rgfl (b)+1}{p_c})
\qnum{\rgfl (b)+1}{p_c}\paar{\iota (\tau (c))}{\tau (c)}
\tau (c)^{\rgfl (b)}\\
&\quad =\paar{\iota (\tau (c))}{\tau (c)}\big(\lambda (b)\tau (\rgf{b})
\tau (a)-\lambda (b)\chi (\rgf{b},a)\tau (a)\tau (\rgf{b})\\
&\qquad -\mu (b)/\qfact{\rgfl (b)}{p_c}\tau (c)^{\rgfl (b)}\big).
\end{align*}
If $\rgfl (b)=1$ then $\rgf{b}=\rgf{c}$, $\tau (\rgf{b})\tau (a)-
\chi (\rgf{b},a)\tau (a)\tau (\rgf{b})=\tau (c)$, and $\mu (b)=\lambda (b)$.
Otherwise $\mu (b)=\lambda (b)\mu (\rgf{b})$ and we can use induction
hypothesis (c). In both cases one obtains $\paar{\iota (\tau (c))}{U}=0$.
Thus by (f) one gets $\paar{\iota (\tau (d))}{U}=0$ for all
$d\in \Nbar _2(T)$ with $c\le _Qd$. Hence it suffices to show that
$\paar{\iota (\suplet{u'})}{U}=0$ for all Lyndon words $u'$ with
$\deg (\suplet{u'})=\deg (U)$ and that $\paar{\iota (\tau (a_1)
\tau (a_2)\cdots \tau (a_m))}{U}=0$ for $a_i\in \Nbar _2(T)$,
$a_1\le _Qa_2\le _Q\cdots \le _Qa_m<_Qc$, $m\ge 2$,
$\sum _i\StBr _T(a_i)=(\rgfl (b)+1)\StBr _T(c)$.
In the second case set $\StBr _T(a_i):=(r_i,s_i)$. Then $a_i<_Qc$
implies $r_i<Q(c)s_i$ and hence $\sum _ir_i<Q(c)\sum _is_i$ which is a
contradiction to $\sum _i\StBr _T(a_i)=(\rgfl (b)+1)\StBr _T(c)$.
In the first case $u'\not=\gamma (d)$ for all $d$ since the entries of
$\deg (\suplet{u'})$ are not relatively prime. Thus by (d) one can reduce
to the case where $u'=v'w'$ is the Shirshow decomposition of $u'$,
$v'=\gamma (a')$, $w'=\gamma (b')$, and $a'=\lgf{\lgf{b'}}$. Since
$(\rgfl (b)+1)\StBr _T(c)=\deg (U)=\deg (\suplet{v'}\suplet{w'})=
(\rgfl (b')+1)\StBr _T(\lgf{b'})$ and both the entries of $\StBr _T(c)$
and those of $\StBr _T(\lgf{b'})$ are relatively prime one gets
from Lemma \ref{l-fbtreeord}(iii) the equation
$c=\lgf{b'}$. Again by Lemma \ref{l-fbtreeord}(iii) this yields
$a=\lgf{c}=\lgf{\lgf{b'}}=a'$ and hence $b=b'$ and $u'=u$.

Note that $\paar{\iota (\tau (b))}{U}=0$ and hence it remains to check
that the relation $\paar{\iota (\tau (b)\tau (a))}{U}=0$ holds.
Using $\tau (b)=\tau (\rgf{b})\tau (c)-\chi (\rgf{b},c)\tau (c)\tau (\rgf{b})$
this implies that it remains to check the equation
$\paar{\iota (\tau (\rgf{b})\tau (\lch{c}))}{U}=0$.
Set $c_1:=\lch{c}$ and $c_{i+1}:=\rch{c_i}$ for $1\le i<\rchl (c_1)$.
Further, set $b_1:=\rgf{b}$ and $b_{i+1}:=\rgf{b_i}$
for $1\le i<\rgfl (b)$. Then using (c) for $(T,V,n-1)$,
the definition of $\tau $, and Lemma \ref{l-tree2Lyndon}(iii) one can show
by induction over $i$ that for all $i\le \min \{\rgfl (b),\rchl (c_1)\}$ the
assertions $\paar{\iota (\tau (b)\tau (a))}{U}=0$ and
$\paar{\iota (\tau (b_i)\tau (c_i))}{U}=0$ are equivalent.
Now we use that $(T,V,n)$ is admissible. More exactly, we have either
$\rchl (c_1)\le \rgfl (b)$ or $\rgfl (b)\le 2$. In the first case one has
$\tau (c_i)=0$ for $i=\rchl (c_1)$. Otherwise (h)
gives $\paar{\iota (\tau (c_i))}{U}=0$ for $i=\rgfl (b)-1$ where
$c_0=a$. Therefore
in both cases we get $\paar{\iota (\tau (b)\tau (a))}{U}=0$ and hence
$U=0$.   

To (b). By (c), (d), (e), and (f) it suffices to show that the equations
$\paar{\iota (\tau (a_1)\tau (a_2)\cdots \tau (a_m))}{\tau (a)}=0$
hold for all $a_1\le _Qa_2\le _Q\cdots
\le _Qa_m\le _Qa$ with $\sum _i\StBr _T(a_i)=\StBr _T(a)$.
Set $(r_i,s_i):=\StBr _T(a_i)$. Then $a_i\le _Qa$ implies $r_i\le Q(a)s_i$
for all $i$ with equality if and only if $a_i=a$. Hence $\sum _i\StBr _T(a_i)
=\StBr _T(a)$ implies $m=1$ and $a_1=a$. Note that $(T,V,n)$ is admissible
and hence $\lambda (a)=0$. Thus the equation
$\spaar{a}{a}=0$ follows from (g).

To (a). Let $a\in \Nbar _2(T)$ and set $u:=\gamma (a)$.
Suppose that $\suplet{u}^{h}$, $h\in \ndN $, is a linear combination of
elements $\suplet{u_1}\suplet{u_2}\cdots \suplet{u_m}$ with
$u<u_i$ for $1\le i\le m$. By Corollary \ref{f-S(u)basis} and the induction
hypothesis one can assume that $u_i=\gamma (a_i)$ for some $a_i$ with
$a<_Qa_i$. Set $(r_i,s_i):=\StBr _T(a_i)$ for all $i$ and
$(r,s):=\StBr _T(a)$. Then Lemma \ref{l-tree2Lyndon} implies that
$r_i>Q(a)s_i$ for all $i$ and hence $\left(\sum _ir_i\right)>
\left(\sum _is_i\right)Q(a)$. On the other hand, since
$\tau (a_1)\tau (a_2)\cdots \tau (a_m)$ must have the same degree
as $\tau (a)^h$ it follows that
$(hr,hs)=(\sum _ir_i,\sum _is_i)$ and hence $\tau (a)^h$ has to be zero.
In particular, $0=\paar{\iota (\tau (a))}{\tau (a)^h}=\qnum{h}{p_a}
\paar{\iota (\tau (a))}{\tau (a)}\tau (a)^{h-1}$ by (e) and (f).
If $n=1$ then $\paar{\iota (\tau (a))}{\tau (a)}=1$. Otherwise
since $(T,V,n)$ is admissible one has $\lambda (a)\not=0$. In this case
(g) gives $\paar{\iota (\tau (a))}{\tau (a)}\not=0$ and hence the first
part of (a) is proven. To show that $\tau (a)^{\mathrm{ord}\,p_a}=0$ for
$n\ge |\StBr _T|(a)\cdot \mathrm{ord}\,p_a$ by (c) and (d) it suffices
to check that $\paar{\iota (\tau (a_1)\tau (a_2)\cdots
\tau (a_m))}{\tau (a)^{\mathrm{ord}\,p_a}}=0$ whenever $a_1\le _Qa_2\le _Q
\cdots \le _Qa_m<_Qa$, $a_i\in \Nbar _2(T)$ for all $i$, and
$\sum _i\StBr _T(a_i)=\StBr _T(a)\cdot \mathrm{ord}\,p_a$.
However as arqued at the beginning of the proof of (a)
such a choice of $a_i$ is not possible.

It remains to prove the induction step for (e)--(h) ($n\to n+1$) under
the hypothesis (a)--(h) and admissibility of $(T,V,n+1)$.

\begin{lemma}\label{l-qkomminF}
For $a\in N(T)$ with $\rgf{a}\in N(T)$, $|\StBr _T|(a)=n+1$ set
$u:=\gamma (a)$ and let $u=vw$ be the Shirshow decomposition of $u$.
Let $\rho \in \cB (V)^+\cap \cF (\gamma (\rgf{\rgf{a}}))_{|w|-1}$ be
a homogeneous element with respect to the $\ndZ ^2$-grading.
Then $g(\suplet{w})^{-1}g(\rho )\iota(\rho )\iota (\suplet{v})
-\chi (\suplet{w},\suplet{v})\iota (\suplet{v})g(\suplet{w})^{-1}g(\rho )
\iota (\rho )\in \iota (\cF (w)_n)$.
\end{lemma}

\begin{bew}[ of the Lemma]
First note that
\begin{gather*}
\begin{aligned}
g(\suplet{w})^{-1}g(\rho )\iota(\rho )\iota (\suplet{v})
-\chi (\suplet{w},\suplet{v})\iota (\suplet{v})g(\suplet{w})^{-1}g(\rho )
\iota (\rho )&\\
=g(\suplet{w})^{-1}g(\rho )\iota (\rho \suplet{v}-
\chi (\rho ,\suplet{v})&\suplet{v}\rho ),
\end{aligned}\\
\tag{$*$}
\begin{aligned}
\rho _1\rho _2\suplet{v}-&\chi (\rho _1\rho _2,\suplet{v})
\suplet{v}\rho _1\rho _2\\
=&\rho _1(\rho _2\suplet{v}-\chi (\rho _2,\suplet{v})\suplet{v}\rho _2)
+\chi (\rho _2,\suplet{v})
(\rho _1\suplet{v}-\chi (\rho _1,\suplet{v})\suplet{v}\rho _1)\rho _2
\end{aligned}
\end{gather*}
for $\ndZ ^2$-homogeneous elements $\rho _i\in \cB (V)$, $i=1,2$.
As $\rho \in S(v)_<$ Lemma \ref{l-freerels} implies that $\rho \suplet{v}
-\chi (\rho ,\suplet{v})\suplet{v}\rho \in S(v)_<$. Further,
$\rho \suplet{v},\suplet{v}\rho \in \cB (V)_m$ where
$m=\mathrm{totdeg}(\rho )+|v|\le n$. By Lemma
\ref{l-fbtreeord}(ix) there exists no $b\in \Nbar _2(T)$ such that
$|\StBr _T|(b)\le n$ and $\lgf{a}<_Qb<_Q\rgf{a}$. Thus
$S(v)_<\cap \cB (V)_m=S(w)\cap \cB (V)_m$.
By ($*$) and since $\rho \in S(\gamma (\rgf{\rgf{a}}))^+\subset S(w)_<^+$
it suffices to show that
\begin{align*}\tag{$**$}
\rho \suplet{v}-\chi (\rho ,\suplet{v})\suplet{v}\rho \in
S(w)S(w)_<^++S(w)\suplet{w}^2
\end{align*}
for $\rho =\tau (\rgf{\rgf{a}})^2$ and for
$\rho =\tau (b)$, $\rgf{\rgf{a}}<_Qb$.
By Equation (\ref{eq-S(u)zerl}) 
relation ($**$) is obviously true if $|v|\ge |w|$. Otherwise
$a$ is the left child of $\rgf{a}$ and $\lgf{\rgf{a}}=\lgf{a}$.
If $\rho =\tau (\rgf{\rgf{a}})^2$ then
$\mathrm{totdeg}(\rho \suplet{v})=|v|+2|\StBr _T|(\rgf{\rgf{a}})
=|w|+|\StBr _T|(\rgf{\rgf{a}})>|w|$
and hence ($**$) holds. On the other hand
$\deg (\suplet{v}\tau (\rgf{\rgf{a}}))=\deg (\suplet{w})$
and Lemma \ref{l-fbtreeord}(iii) imply that
$\deg (\tau (b))\not=\deg (\suplet{w})-\deg (\suplet{v})$
for $b\not=\rgf{\rgf{a}}$. Thus if $\rho =\tau (b)$,
$\rgf{\rgf{a}}<_Qb$ then again ($**$) is valid. 
\end{bew}

To (e). Suppose that $|\StBr _T|(a)=n+1$.
One has $a<_Q\rgf{a}$ and if $\lgf{a}\in N(T)$ then
$\rgf{a}\le _Q\rgf{\lgf{a}}$ by Lemma \ref{l-fbtreeord}(vii),(vi).
Therefore the induction hypothesis (e) for $v$ and $w$ and Lemma
\ref{l-qkomminF} give
\begin{align*}
\copr (\iota (\suplet{v}&))-\iota (\suplet{v})\ot 1-g(\lgf{a})^{-1}\ot
\iota (\suplet{v})\\
-&\lambda (\lgf{a})g(\lgf{\lgf{a}})^{-1}
\iota (\tau (\rgf{\lgf{a}}))\ot \iota (\tau (\lgf{\lgf{a}}))
\in \iota (\cF (w)_{|v|-1})\ot \cB (V^*),\\
\copr (\iota (\suplet{u}&))=
\copr (\iota (\suplet{w}))\copr (\iota (\suplet{v}))
-\chi (\rgf{a},\lgf{a})\copr (\iota (\suplet{v}))
\copr (\iota (\suplet{w}))\\
=&\iota (\suplet{u})\ot 1
+(\chi (\lgf{a},\rgf{a})^{-1}-\chi (\rgf{a},\lgf{a}))
g(\lgf{a})^{-1}\iota (\suplet{w})\ot \iota (\suplet{v})\\
&+g(a)^{-1}\ot \iota (\suplet{u})\\
&+\lambda (\lgf{a})g(\suplet{w}\tau (\lgf{\lgf{a}}))^{-1}
\iota (\tau (\rgf{\lgf{a}}))\ot \iota \big(
\suplet{w}\tau (\lgf{\lgf{a}})-\chi (\rgf{a},\lgf{\lgf{a}})
\tau (\lgf{\lgf{a}})\suplet{w}\big)\\
&+\lambda (\rgf{a})g(\lgf{\rgf{a}})^{-1}\iota \big(
\tau (\rgf{\rgf{a}})\suplet{v}-\chi (\rgf{\rgf{a}},\lgf{a})
\suplet{v}\tau (\rgf{\rgf{a}})\big)\ot \iota (\tau (\lgf{\rgf{a}}))
\end{align*}
up to terms in $\iota (\cF (w)_n)\ot \cB (V^*)$.
Thus if $a=\roo (T)$
then $\lambda (\lgf{a})=\lambda (\rgf{a})=0$ and one gets (\ref{eq-copr[u]}).
If $a$ is the left child of its parent then
either $\lgf{a}=\lgh $ or $\rgf{a}<_Q\rgf{\lgf{a}}$ by Lemma
\ref{l-fbtreeord}(vi),(vii). Further,
$\lgf{\rgf{a}}=\lgf{a}$ and
$\tau (\rgf{\rgf{a}})\suplet{v}-\chi (\rgf{\rgf{a}},\lgf{a})
\suplet{v}\tau (\rgf{\rgf{a}})=\suplet{w}$ and again we are done.
Finally if $a$ is the right child of its parent then
either $\rgf{a}=\rgh $ and $\lambda (\rgf{a})=0$ or
$|v|>|w|$ and
$\tau (\rgf{\rgf{a}})\suplet{v}-\chi (\rgf{\rgf{a}},\lgf{a})
\suplet{v}\tau (\rgf{\rgf{a}})\in \cF (w)_n$ by (\ref{eq-S(u)zerl}).
Moreover $\rgf{\lgf{a}}=\rgf{a}$ and $\suplet{w}\tau (\lgf{\lgf{a}})
-\chi (\rgf{a},\lgf{\lgf{a}})\tau (\lgf{\lgf{a}})\suplet{w}=
\suplet{v}$. Thus (\ref{eq-copr[u]}) holds in this case as well.

To (f).
If $|\StBr _T|(a)\le n$ and $|\StBr _T|(b)\le n$ then we are done by induction
hypothesis. If $n+1=|\StBr _T|(b)>|\StBr _T|(a)$ then
$\paar{\iota (\tau (b))}{\tau (a)}=0$ clearly holds.
Suppose now that $|\StBr _T|(a)=n+1$. If $n=0$ then $a=\lgh $ and $b=\rgh $
and hence (f) holds by definition of $\paar{\cdot }{\cdot }$.
Otherwise $\rgf{a}\le _Qb$ by Lemma \ref{l-fbtreeord}(vi) and hence (e) and
(f) imply that 
\begin{align*}
\paar{\iota (\tau (b))}{\tau (a)}=&\paar{\iota (\tau (b))}{\tau (\rgf{a})
\tau (\lgf{a})-\chi (\rgf{a},\lgf{a})\tau (\lgf{a})\tau (\rgf{a})}\\
=&\paar{\iota (\tau (b))}{\tau (\rgf{a})}\lgf{a}
-\chi (\rgf{a},\lgf{a})\chi (b,\lgf{a})^{-1}\lgf{a}
\paar{\iota (\tau (b))}{\tau (\rgf{a})}.
\end{align*}
By (f) the last expression vanishes if $\rgf{a}\le _Qb$.

To (g). By (f) and since $\tau (a)=\tau (\rgf{a})\tau (\lgf{a})
-\chi (\rgf{a},\lgf{a})\tau (\lgf{a})\tau (\rgf{a})$ the second equation
of (g) follows immediately from the first one.

If $a=\roo (T)$ then one gets $\paar{\iota (\tau (\lgf{a}))}{\tau (a)}=
\paar{y_2}{x_1x_2-\chi (x_1,x_2)x_2x_1}=(\chi (x_2,x_1)^{-1}-\chi (x_1,x_2))
x_1=\lambda (a)\tau (\rgf{a})$.

If $a=\rch{\lgf{a}}$ then $\rgf{\lgf{a}}=\rgf{a}$ and $|\StBr _T|(\lgf{a})
>|\StBr _T|(\rgf{a})$ and hence (e) and (f) gives
\begin{align*}
\paar{\iota (\tau (\lgf{a}))}{\tau (a)}=&
\paar{\iota (\tau (\lgf{a}))}{\tau (\rgf{a})\tau (\lgf{a})-\chi (\rgf{a},
\lgf{a})\tau (\lgf{a})\tau (\rgf{a})}\\
=&\chi (\lgf{a},\rgf{a})^{-1}
\tau (\rgf{a})\paar{\iota (\tau (\lgf{a}))}{\tau (\lgf{a})}\\
&+\lambda (\lgf{a})\paar{\iota (\tau (\rgf{a}))}{\tau (\rgf{a})}
\paar{\iota (\tau (\lgf{\lgf{a}}))}{\tau (\lgf{a})}\\
&-\chi (\rgf{a},\lgf{a})\paar{\iota (\tau (\lgf{a}))}{\tau (\lgf{a})}
\tau (\rgf{a}).
\end{align*}
Since $\lgf{a}\in N_2(T)$ in the second summand of the last expression
one can use (g) for $\lgf{a}$ and equation $\rgf{\lgf{a}}=\rgf{a}$.
One gets $\paar{\iota (\tau (\lgf{a}))}{\tau (a)}=(\lambda (\lgf{a})+
\chi (\lgf{a},\rgf{a})^{-1}-\chi (\rgf{a},\lgf{a}))
\paar{\iota (\tau (\lgf{a}))}{\tau (\lgf{a})}\tau (\rgf{a})
=\lambda (a)\paar{\iota (\tau (\lgf{a}))}{\tau (\lgf{a})}\tau (\rgf{a})$.

Finally, if $a=\lch{\rgf{a}}$ then $\lgf{\rgf{a}}=\lgf{a}$. Moreover
Lemma \ref{l-fbtreeord}(vii) implies that $\lgf{a}=\lgh $ or
$\rgf{a}<_Q\rgf{\lgf{\rgf{a}}}=\rgf{\lgf{a}}$. Thus using (e), (f), and (g)
computations similar to the previous case lead to the desired assertion.

To (h). We need the following lemma.

\begin{lemma}\label{l-bigcontr}
Suppose that $(T,V,n+1)$ is admissible. Let $d\in N_2(T)$ such that
$c:=\rgf{d}\in N(T)$ and $d=\lch{c}$. Set $a:=\lgf{c}=\lgf{d}$ and
$f:=\rgf{c}$ (see Fig.~2).\\
(i) If $|\StBr _T|(d)+|\StBr _T|(f)\le n$ then in $\cB (V)$ one has
\begin{align*}
\tau (f)\tau (d)-\chi (f,d)\tau (d)\tau (f)=
\frac{\lambda (d)}{\qnum{2}{p_c}}\tau (c)^2.
\end{align*}
(ii) If $|\StBr _T|(d)+|\StBr _T|(f)(=|\StBr _T|(a)+|\StBr _T|(\rch{c}))
\le n+1$ then
\begin{align*}
\paar{\iota (\tau (a))}{\tau (\rch{c})}=\frac{\lambda (c)\lambda (\rch{c})}{
\qnum{2}{p_f}}\paar{\iota (\tau (a))}{\tau (a)}\tau (f)^2.
\end{align*}
(iii) If $\rch{c}\in N_2(T)$ and $|\StBr _T|(\lch{\rch{c}})+
|\StBr _T|(a)\le n+1$ then
\begin{align*}
\paar{\iota (\tau (d))}{\tau (\lch{\rch{c}})}=\frac{\lambda (\rch{c})
\lambda (\lch{\rch{c}})\lambda (c)
\paar{\iota (\tau (d))}{\tau (d)}}{\qnum{2}{p_c}\qnum{2}{p_f}}\tau (f)^2.
\end{align*}
\end{lemma}

\begin{bew}[ of the Lemma]
To (i). Note that $\lgf{c}=\lgf{d}$ by assumption. Further,
we may use the induction hypothesis (c). One computes
{\allowdisplaybreaks
\begin{align*}
\tau (f)\tau (d)&-\chi (f,d)\tau (d)\tau (f)\\
=&\tau (f)(\tau (c)\tau (a)-\chi (c,a)\tau (a)\tau (c))
-\chi (f,d)\tau (d)\tau (f)\\
=&(\tau (\rch{c})+\chi (f,c)\tau (c)\tau (f))\tau (a)\\
&-\chi (c,a)(\tau (c)+\chi (f,a)\tau (a)\tau (f))\tau (c)
-\chi (f,d)\tau (d)\tau (f)\\
=&\chi (\rch{c},a)\tau (a)\tau (\rch{c})+\lambda (\rch{c})/\qnum{2}{p_c}
\tau (c)^2\\
&+\chi (f,c)\tau (c)(\tau (c)+\chi (f,a)\tau (a)\tau (f))
-\chi (c,a)\tau (c)^2\\
&-\chi (\rch{c},a)\tau (a)(\tau (\rch{c})+\chi (f,c)\tau (c)\tau (f))
-\chi (f,d)\tau (d)\tau (f)\\
=&\lambda (\rch{c})/\qnum{2}{p_c}\tau (c)^2+\chi (f,c)\tau (c)^2\\
&+\chi (f,d)(\tau (d)+\chi (c,a)\tau (a)\tau (c))\tau (f)
-\chi (c,a)\tau (c)^2\\
&-\chi (\rch{c},a)\chi (f,c)\tau (a)\tau (c)\tau (f)
-\chi (f,d)\tau (d)\tau (f)\\
=&(\lambda (\rch{c})+\chi (f,c)+\chi (a,c)^{-1}
-\chi (c,a)-\chi (c,f)^{-1})/\qnum{2}{p_c}\tau (c)^2.
\end{align*}
}
Then the defining recursion formulas for $\lambda (\rch{c})$ and $\lambda (d)$
give (i).

To (ii). Note that $\qnum{2}{p_f}\not=0$ by admissibility of
$(T,V,n+1)$. Further, if $\tau (\rch{c})=0$ then $\lambda (\rch{c})=0$ by (a).
Thus in this case we are done. Assume now that $c=\roo $
or $\rch{a}=c$ (i.\,e.~$\rgf{a}=f$). Then one has
$\paar{\iota (\tau (a))}{\tau (f)}=0$. Using (e) and (f)
one gets
\begin{align*}
&\spaar{a}{\rch{c}}=
\paar{\iota (\tau (a))}{\tau (f)\tau (c)-\chi (f,c)\tau (c)\tau (f)}\\
&\quad =\chi (a,f)^{-1}\tau (f)\spaar{a}{c}+\lambda (a)\spaar{f}{f}
\spaar{\lgf{a}}{c}\\
&\qquad -\chi (f,c)\spaar{a}{c}\tau (f).
\end{align*}
If one starts with $\lch{a}$ instead of $d$ in Lemma \ref{l-bigcontr} then (ii)
gives a formula for the second summand in the last expression.
Further, (g) can be used to compute $\spaar{a}{c}$ and $\spaar{a}{a}$.
One obtains
{\allowdisplaybreaks
\begin{align*}
&\spaar{a}{\rch{c}}=
\lambda (c)\spaar{a}{a}(\chi (a,f)^{-1}-\chi (f,c))\tau (f)^2\\
&\qquad +\lambda (a)\spaar{f}{f}\frac{\lambda (c)\lambda (a)}{\qnum{2}{p_f}}
\spaar{\lgf{a}}{\lgf{a}}\tau (f)^2\\
&\quad =\frac{\lambda (c)\spaar{a}{a}}{\qnum{2}{p_f}}
(\chi (a,f)^{-1}{+}\chi (c,f)^{-1}{-}\chi (f,c){-}\chi (f,a){+}\lambda (a))
\tau (f)^2.
\end{align*}
}
Note that these computations make sense also in the case when $a=\lgh $.
Thus the recursion formulas for $\lambda (c)$ and $\lambda (\rch{c})$ give
(ii) in this case. The proof of (ii) in the remaining case (when $\lch{f}=c$,
i.\,e.~$\lgf{c}=\lgf{f}$,) is obtained similarly.

The proof of (iii) is by far the most complicated one.
We give only a scetch of it.
Set $b:=\lch{\rch{c}}$. Using (e), (f), (g), and Lemma \ref{l-bigcontr}(ii)
one obtains
\begin{align*}
\spaar{a}{b}=&\frac{\lambda (c)\spaar{a}{a}}{\qnum{2}{p_f}}\bigg(
(\lambda (\rch{c})-\chi (\rch{c},c)\qnum{2}{p_f})\tau (f)\tau (\rch{c})\\
&+(\chi (f,c)\lambda (\rch{c})+\chi (a,\rch{c})^{-1}\qnum{2}{p_f})
\tau (\rch{c})\tau (f)\bigg),\\
\paar{\iota (\tau (c))}{\tau (f)\tau (\rch{c})}=&
\frac{\lambda (\rch{c})\spaar{c}{c}}{\qnum{2}{p_f}}(\lambda (c)+
\chi (c,f)^{-1}\qnum{2}{p_f})\tau (f)^2,\\
\paar{\iota (\tau (c))}{\tau (\rch{c})\tau (f)}=&
\lambda (\rch{c})\spaar{c}{c}\tau (f)^2,\\
\paar{\iota (\tau (a)\tau (c))}{\tau (b)}=&
\frac{\lambda (\rch{c})\lambda (b)\lambda (c)}{\qnum{2}{p_f}}
\spaar{c}{c}\spaar{a}{a}\tau (f)^2.
\end{align*}
By (c) and (h) one obtains $\nu (\rch{c})=0$. Now insert the second, third,
and fourth equation into the first one and replace the summand
$\lambda (\rch{c})\lambda (c)\qnum{2}{p_f}^{-1}$ of the product
$(\lambda (\rch{c})-\chi (\rch{c},c)\qnum{2}{p_f})(\lambda (c)+
\chi (c,f)^{-1}\qnum{2}{p_f})$ by the expression
$\lambda (\rch{c})\lambda (c)\qnum{2}{p_c}^{-1}-\chi (a,\rch{c})^{-1}
+\chi (\rch{c},a)$. Then use the recursion formulas for $\lambda (d)$,
$\lambda (b)$, and $\lambda (\rch{c})$ to obtain (iii).
\end{bew}

To show the first equation of (h) one can use (e) and (f) to obtain
\begin{align*}
&\paarb{\iota (\tau (a))}{\tau (b)\tau (a)-\chi (b,a)\tau (a)\tau (b)
-\frac{\lambda (b)}{\qnum{2}{p_c}}\tau (c)^2}=\\
&\quad =\spaar{a}{b}\tau (a)+\chi (a,b)^{-1}\tau (b)\spaar{a}{a}\\
&\qquad 
-\chi (b,a)\spaar{a}{a}\tau (b)-\chi (f,a)^2\tau (a)\spaar{a}{b}\\
&\qquad -\frac{\lambda (b)}{\qnum{2}{p_c}}(\spaar{a}{c}\tau (c)
+\chi (a,c)^{-1}\tau (c)\spaar{a}{c}).
\end{align*}
After inserting the formula in Lemma \ref{l-bigcontr}(ii) the latter
expression becomes the sum of
$\spaar{a}{a}\nu (b)\tau (b)$ and a multiple of $\tau (c)\tau (f)$. Apply
$\paar{\iota (\tau (f))}{\cdot }$ to this equation.
Since one has $\spaar{f}{b}=0$ and $\paar{\iota (\tau (f))}{\tau (c)\tau (f)}=
\chi (f,c)^{-1}\tau (c)\spaar{f}{f}$ Lemma \ref{l-bigcontr}(i) gives that
the coefficient of $\tau (c)\tau (f)$ is zero.

In order to prove the second equation of (h) note that
$|\StBr _T|(d)=|\StBr _T|(c)+
|\StBr _T|(a)<2|\StBr _T|(c)<|\StBr_T|(b)$. Thus (e) gives that
\begin{align*}
\copr (\iota (\tau (d)))-\iota (\tau (d))\ot 1-g(d)^{-1}\ot
\iota (\tau (d))-\lambda (d)g(a)^{-1}\iota (\tau (c))\ot \iota (\tau (a))&\\
\in \Big(\iota \big(S(\gamma (\rch{c}))^+ +\tau (c)S(\gamma (\rch{c}))^+\big)
\# kG\Big)\ot \cB (V^*)&. 
\end{align*}
By this fact and assumption (f) one obtains that
\begin{align*}
&\paarb{\iota (\tau (d))}{\tau (b)\tau (a)-\chi (b,a)\tau (a)\tau (b)
-\frac{\lambda (b)\lambda (\rch{c})}{\qfact{3}{p_c}}\tau (c)^3}=\\
&\quad =\spaar{d}{b}\tau (a)
+\lambda (d)\paar{g(a)^{-1}\iota (\tau (c))}{\tau (b)}\spaar{a}{a}\\
&\qquad -\chi (b,a)\chi (d,a)^{-1}\tau (a)\spaar{d}{b}\\
&\qquad -\frac{\lambda (b)\lambda (\rch{c})}{\qfact{3}{p_c}}\Big(
\chi (d,c)^{-1}\tau (c)\paar{\iota (\tau (d))}{\tau (c)^2}\\
&\qquad +\lambda (d)\spaar{c}{c}\paar{\iota (\tau (a))}{\tau (c)^2}\Big).
\end{align*}
Using Lemma \ref{l-bigcontr}(iii) and arguments as in the previous case
one obtains the required result.
\end{bew}

\begin{bew}[ of Theorem \ref{t-class}]
The proof bases on Corollary \ref{f-findeg} and consists of a case by case
checking of admissibility of $(T,V,n)$ for all $n\in \ndN $.

Suppose that $a\in N_2(T)$. Let $b,c\in N_0(T)$ be the unique nodes
such that $\rgf{b}=\lgf{c}=a$. Then by (\ref{eq-lambda}) one obtains
that
\begin{align*}
\lambda (b)=&\lambda (a)+(\chi (\lgf{a},a)^{-1}-\chi (a,\lgf{b}))
\qnum{\lgfl (b)}{p_a},\\
\lambda (c)=&\lambda (a)+(\chi (a,\rgf{a})^{-1}-\chi (\rgf{c},a))
\qnum{\rgfl (c)}{p_a}.
\end{align*}
Thus one gets
\begin{align}
\lambda (b)=\lambda (c) \Leftrightarrow &
\chi (\lgf{a},a)^{-1}\qnum{\lgfl (b)}{p_a}
+\chi (\rgf{c},a)\qnum{\rgfl (c)}{p_a}= \notag \\
& =\chi (a,\lgf{b})\qnum{\lgfl (b)}{p_a}
+\chi (a,\rgf{a})^{-1}\qnum{\rgfl (c)}{p_a} \notag \\
\Leftrightarrow &\qnum{\lgfl (b)+\rgfl (c)}{p_a}(\chi (\rgf{c},a)-
\chi (a,\lgf{b}))=0 \notag \\
\label{eq-fin1}
\Leftrightarrow & \qnum{\lgfl (b)+\rgfl (c)}{p_a}(p_{\rgf{a}}
-p_a^{\lgfl (b)-\rgfl (c)}p_{\lgf{a}})=0.
\end{align}
Further, if $a\in N(T)$ and $\rgf{a}=\rgh $ then
\begin{align}\label{eq-fin2}
\lambda (a)=(q_{11}q_{12}q_{22})^{-1}
(p_{\roo }-q_{22}^{-1}/q_{11}^{-(\lgfl (a)-2)})\qnum{\lgfl (a)}{q_{11}^{-1}},
\end{align}
and if $a\in N(T)$ and $\lgf{a}=\lgh $ then
\begin{align}\label{eq-fin3}
\lambda (a)=(q_{11}q_{12}q_{22})^{-1}(p_{\roo }-q_{11}^{-1}/
q_{22}^{-(\rgfl (a)-2)})\qnum{\rgfl (a)}{q_{22}^{-1}}.
\end{align}
Equations (\ref{eq-fin1}), (\ref{eq-fin2}), and (\ref{eq-fin3}) give an
effective method to check the first condition of
admissibility of $(T,V,n)$ for all $n$. In fact, the equivalence between
$\lambda (a)=0$ and $a\in N_0(T)$ holds for all $a\in N(T)$
if and only if
\begin{align*}
\rchl (\roo )=&\min \{m\in \ndN \,|\,\qnum{m}{q_{11}^{-1}}(q_{11}^{1-m}
p_{\roo }-q_{11}^{-1}q_{22}^{-1})=0\},\\
\lchl (\roo )=&\min \{m\in \ndN \,|\,\qnum{m}{q_{22}^{-1}}(q_{22}^{1-m}
p_{\roo }-q_{22}^{-1}q_{11}^{-1})=0\},\\
\rchl (\lch{a})=&\min \{m\in \ndN \,|\,\qnum{m+\lchl (\rch{a})}{p_a}
(p_{\rgf{a}}p_a^{\lchl (\rch{a})}-p_{\lgf{a}}p_a^m)=0\}
\,\forall a\in N_2(T).
\end{align*}

Let $a_1,a_2,\ldots $ denote the elements of $N_2(T)$ such that
$a_i<_Qa_j$ for $i<j$ and set $p_i:=p_{a_i}$. In what follows we give
all values $p_i$ and all $\lambda (\cdot )$ which are relevant for the
computation of the necessary $\nu (a)$. However if $\rgf{a}=\rgh $ then
one has the closed formula
$\lambda (a)=q_{21}^{-1}\qnum{\lgfl (a)}{q_{11}^{-1}}
(1-q_{11}^{\lgfl (a)-1}q_{12}q_{21})$
and we will omit to give $\lambda (a)$ more explicitly.\\
(T1) is trivial.\\
(T2) $\lambda (\lch{\roo })=(1+q_{22}^{-1})(q_{21}^{-1}-q_{12}q_{22})$,
$p_1=(q_{11}q_{12}q_{21}q_{22})^{-1}$.\\
(T3) In the case $q_{12}q_{21}=q_{11}^2$ one has
$p_1=q_{11}/q_{22}$, $p_2=q_{22}^{-1}$.\\
If $\qnum{3}{q_{11}^{-1}}=0$ and $q_{12}q_{21}q_{22}=1$ then $p_1=q_{11}^{-1}$,
$p_2=q_{22}/q_{11}$.\\
In the third case one gets $p_1=q_{11}$, $p_2=-1$.\\
(T4) First case: $p_1=-1$, $p_2=q_0^3$, $p_3=-1$.\\
Second case: $p_1=-1$, $p_2=-q_{12}q_{21}$, $p_3=-1$.\\
(T5) First case: $p_1=-(q_{12}q_{21})^3$, $p_2=-1$, $p_3=-(q_{12}q_{21})^2$.\\
Second case: $p_1=q_0^5$, $p_2=-1$, $p_3=-q_0^2$.\\
(T6) $p_1=-1$, $p_2=-q_{11}^{-2}$, $p_3=-1$, $p_4=-q_{11}^{-3}$.\\
(T7) First case: $p_1=-q_{11}^2$, $p_2=-q_{11}^2$, $p_3=-1$.\\
Second case: $p_1=-(q_{12}q_{21})^2$, $p_2=(q_{12}q_{21})^4$, $p_3=-1$.\\
(T8) First case: $p_1=q_{11}^{-1}$, $p_2=q_{11}^{-3}$, $p_3=q_{11}^{-1}$,
$p_4=q_{11}^{-3}$.\\
Second case ($q_{11}=-q_{12}q_{21}$): $p_1=-(q_{12}q_{21})^2$,
$p_2=-q_{12}q_{21}$, $p_3=-(q_{12}q_{21})^2$, $p_4=-1$.\\
Third case: ($q_{11}=-(q_{12}q_{21})^2$): $p_1=-q_{12}q_{21}$,
$p_2=-1$, $p_3=(q_{12}q_{21})^2$, $p_4=(q_{12}q_{21})^3$.\\
Fourth case: $p_1=-(q_{12}q_{21})^2$, $p_2=-1$,
$p_3=-(q_{12}q_{21})^3$, $p_4=-1$.\\
(T9) $p_1=-(q_{12}q_{21})^2$, $p_2=-1$, $p_3=(q_{12}q_{21})^3$,
$p_4=-q_{12}q_{21}$.\\
(T10) $p_1=-1$, $p_2=-q_{12}q_{21}$, $p_3=-1$, $p_4=(q_{12}q_{21})^8$,
$p_5=(q_{12}q_{21})^6$, $p_6=-(q_{12}q_{21})^{-1}$.\\
(T11) $p_1=-q_{11}^2$, $p_2=-1$, $p_3=q_{11}^9$, $p_4=-1$,
$p_5=-q_{11}^2$, $p_6=-1$,
$\lambda (a_3)=q_{21}^{-1}(1-q_{11}^{-2}+q_{11}^{-4})(1-q_{11})$,
$\lambda (a_4)=q_{21}^{-1}(1-q_{11}^{-2})(1+q_{11}^{-4})(1+q_{11}^3)$.\\
(T12) $p_1=-1$, $p_2=-q_{11}^{-3}$, $p_3=-q_{11}^4$, $p_4=-q_{11}^{-3}$,
$p_5=-1$, $p_6=-q_{11}^{-5}$.\\
(T13) $p_1=-(q_{12}q_{21})^6$, $p_2=-(q_{12}q_{21})^4$, $p_3=-1$,
$p_4=(q_{12}q_{21})^{-1}$, $p_5=-1$, $p_6=-(q_{12}q_{21})^4$.\\
(T14) $p_1=-q_{11}^3$, $p_2=-q_{11}^4$, $p_3=-q_{11}^3$, $p_4=-1$.\\
(T15) $p_1=-(q_{12}q_{21})^3$, $p_2=-1$, $p_3=(q_{12}q_{21})^{10}$,
$p_4=(q_{12}q_{21})^{11}$, $p_5=(q_{12}q_{21})^{10}$, $p_6=-1$.\\
(T16) First case: $p_1=-q_{11}^3$, $p_2=-1$, $p_3=-q_{11}^4$,
$p_4=-1$, $p_5=-q_{11}^3$, $p_6=-1$.\\
Second case: $p_1=-(q_{12}q_{21})^3$, $p_2=-1$, $p_3=(q_{12}q_{21})^4$,
$p_4=-1$, $p_5=(q_{12}q_{21})^3$, $p_6=-1$.\\
(T17) $p_1=-(q_{12}q_{21})^7$, $p_2=-1$, $p_3=(q_{12}q_{21})^8$,
$p_4=-(q_{12}q_{21})^6$, $p_5=(q_{12}q_{21})^5$, $p_6=-(q_{12}q_{21})^6$.\\
(T18) $p_1=-(q_{12}q_{21})^9$, $p_2=-(q_{12}q_{21})^{-2}$,
$p_3=-(q_{12}q_{21})^9$, $p_4=-1$,
$p_5=(q_{12}q_{21})^{10}$, $p_6=-(q_{12}q_{21})^8$,
$\lambda (a_3)=q_{21}^{-1}(1+q_{12}q_{21})((q_{12}q_{21})^4
+(q_{12}q_{21})^{11})$,
$\lambda (a_5)=q_{21}^{-1}(q^5-q^{-4})$.\\
(T19) $p_1=-q_{11}^2$, $p_2=-1$, $p_3=q_{11}^{-1}$, $p_4=-1$,
$p_5=-q_{11}^2$, $p_6=-1$, $p_7=q_{11}^{-1}$, $p_8=-1$,
$p_9=-q_{11}^2$, $p_{10}=-1$,
$\lambda (a_5)=q_{21}^{-1}(1-q_{11})(1-q_{11}^{-2}+q_{11}^{-4})$,
$\lambda (a_3)=\lambda (a_7)=q_{21}^{-1}(1+q_{11}^{-1})(1-q_{11}^{-2})$,
$\lambda (a_4)=\lambda (a_8)=q_{21}^{-1}(1-q_{11}^{-3})$.\\
(T20) $p_1=-(q_{12}q_{21})^5$, $p_2=(q_{12}q_{21})^7$,
$p_3=-(q_{12}q_{21})^5$, $p_4=-1$, $p_5=(q_{12}q_{21})^6$,
$p_6=-(q_{12}q_{21})^2$.\\
(T21) $p_1=-q_{11}^4$, $p_2=-q_{11}^6$, $p_3=q_{11}$, $p_4=-q_{11}^6$,
$p_5=-q_{11}^4$, $p_6=-1$.\\
(T22) $p_1=-q_{11}^4$, $p_2=-1$, $p_3=q_{11}^{-1}$, $p_4=-1$,
$p_5=-q_{11}^4$, $p_6=-1$, $p_7=q_{11}^{-1}$, $p_8=-1$,
$p_9=-q_{11}^4$, $p_{10}=-1$,
$\lambda (a_5)(=\lambda _9)=q_{21}^{-1}(1+q_{11}^{-1})(1-q_{11}^{-4})$,
$\lambda (a_6)(=\lambda _8=\lambda (a_{10}))=q_{21}^{-1}(1-q_{11}^{-5})$.
\end{bew}

\begin{appendix}
\section{Types of Nichols algebras}
\label{app-types}

In this appendix we collect all full binary trees which appear as the type
of some finite dimensional Nichols algebra studied in the present paper.

\begin{center}
\setlength{\unitlength}{5pt}
\begin{picture}(2,2)
\put(1,1){\circle*{.5}}
\end{picture}
T1
\hfill \hfill
\begin{picture}(4,2)
\put(2,2){\line(-1,-1){2}}
\put(2,2){\line(1,-1){2}}
\put(2,2){\circle*{.5}}
\put(0,0){\circle*{.5}}
\put(4,0){\circle*{.5}}
\end{picture}
\,T2
\hfill \hfill
\begin{picture}(6,4)
\put(0,2){\line(1,1){2}}
\put(2,4){\line(1,-1){4}}
\put(2,0){\line(1,1){2}}
\put(2,4){\circle*{.5}}
\put(0,2){\circle*{.5}}
\put(4,2){\circle*{.5}}
\put(2,0){\circle*{.5}}
\put(6,0){\circle*{.5}}
\end{picture}
\,T3
\hfill \hfill
\begin{picture}(6,4)
\put(3,4){\line(-1,-1){2}}
\put(3,4){\line(1,-1){2}}
\put(1,2){\line(1,-2){1}}
\put(1,2){\line(-1,-2){1}}
\put(5,2){\line(-1,-2){1}}
\put(5,2){\line(1,-2){1}}
\put(3,4){\circle*{.5}}
\put(1,2){\circle*{.5}}
\put(5,2){\circle*{.5}}
\put(0,0){\circle*{.5}}
\put(2,0){\circle*{.5}}
\put(4,0){\circle*{.5}}
\put(6,0){\circle*{.5}}
\end{picture}
\,T4
\end{center}
\begin{center}
\setlength{\unitlength}{5pt}
\begin{picture}(6,6)
\put(0,4){\line(1,1){2}}
\put(1,0){\line(1,2){1}}
\put(2,2){\line(1,1){2}}
\put(2,2){\line(1,-2){1}}
\put(2,6){\line(1,-1){4}}
\put(2,6){\circle*{.5}}
\put(0,4){\circle*{.5}}
\put(4,4){\circle*{.5}}
\put(2,2){\circle*{.5}}
\put(6,2){\circle*{.5}}
\put(1,0){\circle*{.5}}
\put(3,0){\circle*{.5}}
\end{picture}
\,T5
\hfill \hfill
\begin{picture}(8,6)
\put(0,2){\line(1,1){4}}
\put(2,4){\line(1,-2){1}}
\put(4,6){\line(1,-1){4}}
\put(6,4){\line(-1,-2){2}}
\put(5,2){\line(1,-2){1}}
\put(4,6){\circle*{.5}}
\put(2,4){\circle*{.5}}
\put(6,4){\circle*{.5}}
\put(0,2){\circle*{.5}}
\put(3,2){\circle*{.5}}
\put(5,2){\circle*{.5}}
\put(8,2){\circle*{.5}}
\put(4,0){\circle*{.5}}
\put(6,0){\circle*{.5}}
\end{picture}
\,T6
\hfill \hfill
\begin{picture}(8,6)
\put(0,4){\line(1,1){2}}
\put(2,2){\line(1,1){2}}
\put(4,0){\line(1,1){2}}
\put(2,6){\line(1,-1){6}}
\put(2,6){\circle*{.5}}
\put(0,4){\circle*{.5}}
\put(4,4){\circle*{.5}}
\put(2,2){\circle*{.5}}
\put(6,2){\circle*{.5}}
\put(4,0){\circle*{.5}}
\put(8,0){\circle*{.5}}
\end{picture}
\,T7
\hfill \hfill
\begin{picture}(7,6)
\put(0,4){\line(1,1){2}}
\put(1,0){\line(1,2){1}}
\put(2,2){\line(1,1){2}}
\put(2,2){\line(1,-2){1}}
\put(2,6){\line(1,-1){4}}
\put(6,2){\line(-1,-2){1}}
\put(6,2){\line(1,-2){1}}
\put(2,6){\circle*{.5}}
\put(0,4){\circle*{.5}}
\put(4,4){\circle*{.5}}
\put(2,2){\circle*{.5}}
\put(6,2){\circle*{.5}}
\put(1,0){\circle*{.5}}
\put(3,0){\circle*{.5}}
\put(5,0){\circle*{.5}}
\put(7,0){\circle*{.5}}
\end{picture}
\,T8
\end{center}
\begin{center}
\setlength{\unitlength}{5pt}
\begin{picture}(6,8)
\put(0,6){\line(1,1){2}}
\put(1,2){\line(1,2){1}}
\put(2,4){\line(1,1){2}}
\put(2,4){\line(1,-2){1}}
\put(2,8){\line(1,-1){4}}
\put(1,2){\line(-1,-2){1}}
\put(1,2){\line(1,-2){1}}
\put(2,8){\circle*{.5}}
\put(0,6){\circle*{.5}}
\put(4,6){\circle*{.5}}
\put(2,4){\circle*{.5}}
\put(6,4){\circle*{.5}}
\put(1,2){\circle*{.5}}
\put(3,2){\circle*{.5}}
\put(0,0){\circle*{.5}}
\put(2,0){\circle*{.5}}
\end{picture}
\,T9
\hfill \hfill
\begin{picture}(8,8)
\put(1,6){\line(-1,-2){1}}
\put(1,6){\line(1,-2){1}}
\put(3,8){\line(-1,-1){2}}
\put(3,8){\line(1,-1){4}}
\put(5,6){\line(-1,-1){2}}
\put(3,4){\line(-1,-2){1}}
\put(3,4){\line(1,-2){1}}
\put(7,4){\line(-1,-2){1}}
\put(7,4){\line(1,-2){1}}
\put(2,2){\line(-1,-2){1}}
\put(2,2){\line(1,-2){1}}
\put(3,8){\circle*{.5}}
\put(1,6){\circle*{.5}}
\put(5,6){\circle*{.5}}
\put(0,4){\circle*{.5}}
\put(2,4){\circle*{.5}}
\put(3,4){\circle*{.5}}
\put(7,4){\circle*{.5}}
\put(2,2){\circle*{.5}}
\put(4,2){\circle*{.5}}
\put(6,2){\circle*{.5}}
\put(8,2){\circle*{.5}}
\put(1,0){\circle*{.5}}
\put(3,0){\circle*{.5}}
\end{picture}
\,T10
\hfill \hfill
\begin{picture}(8,8)
\put(0,0){\line(1,2){1}}
\put(1,2){\line(1,-2){1}}
\put(1,6){\line(1,1){2}}
\put(3,8){\line(1,-1){4}}
\put(1,2){\line(1,1){4}}
\put(3,4){\line(1,-1){2}}
\put(5,2){\line(-1,-2){1}}
\put(5,2){\line(1,-2){1}}
\put(7,4){\line(-1,-2){1}}
\put(7,4){\line(1,-2){1}}
\put(3,8){\circle*{.5}}
\put(1,6){\circle*{.5}}
\put(5,6){\circle*{.5}}
\put(3,4){\circle*{.5}}
\put(7,4){\circle*{.5}}
\put(1,2){\circle*{.5}}
\put(5,2){\circle*{.5}}
\put(6,2){\circle*{.5}}
\put(8,2){\circle*{.5}}
\put(0,0){\circle*{.5}}
\put(2,0){\circle*{.5}}
\put(4,0){\circle*{.5}}
\put(6,0){\circle*{.5}}
\end{picture}
\,T11
\hfill \hfill
\begin{picture}(8,8)
\put(1,6){\line(-1,-2){1}}
\put(1,6){\line(1,-2){1}}
\put(3,8){\line(-1,-1){2}}
\put(3,8){\line(1,-1){4}}
\put(5,6){\line(-1,-1){2}}
\put(3,4){\line(-1,-2){1}}
\put(3,4){\line(1,-2){1}}
\put(7,4){\line(-1,-1){2}}
\put(7,4){\line(1,-2){1}}
\put(5,2){\line(-1,-2){1}}
\put(5,2){\line(1,-2){1}}
\put(3,8){\circle*{.5}}
\put(1,6){\circle*{.5}}
\put(5,6){\circle*{.5}}
\put(0,4){\circle*{.5}}
\put(2,4){\circle*{.5}}
\put(3,4){\circle*{.5}}
\put(7,4){\circle*{.5}}
\put(2,2){\circle*{.5}}
\put(4,2){\circle*{.5}}
\put(5,2){\circle*{.5}}
\put(8,2){\circle*{.5}}
\put(4,0){\circle*{.5}}
\put(6,0){\circle*{.5}}
\end{picture}
\,T12
\end{center}
\begin{center}
\setlength{\unitlength}{5pt}
\begin{picture}(7,8)
\put(2,8){\line(-1,-1){2}}
\put(2,8){\line(1,-1){4}}
\put(4,6){\line(-3,-2){3}}
\put(1,4){\line(-1,-2){1}}
\put(1,4){\line(1,-2){1}}
\put(6,4){\line(-1,-2){1}}
\put(6,4){\line(1,-2){1}}
\put(2,2){\line(-1,-2){1}}
\put(2,2){\line(1,-2){1}}
\put(5,2){\line(-1,-2){1}}
\put(5,2){\line(1,-2){1}}
\put(2,8){\circle*{.5}}
\put(0,6){\circle*{.5}}
\put(4,6){\circle*{.5}}
\put(1,4){\circle*{.5}}
\put(6,4){\circle*{.5}}
\put(0,2){\circle*{.5}}
\put(2,2){\circle*{.5}}
\put(5,2){\circle*{.5}}
\put(7,2){\circle*{.5}}
\put(1,0){\circle*{.5}}
\put(3,0){\circle*{.5}}
\put(4,0){\circle*{.5}}
\put(6,0){\circle*{.5}}
\end{picture}
\,T13
\hfill \hfill
\begin{picture}(9,8)
\put(2,8){\line(-1,-1){2}}
\put(2,8){\line(1,-1){6}}
\put(4,6){\line(-1,-1){2}}
\put(6,4){\line(-1,-1){2}}
\put(8,2){\line(-1,-2){1}}
\put(8,2){\line(1,-2){1}}
\put(2,8){\circle*{.5}}
\put(0,6){\circle*{.5}}
\put(4,6){\circle*{.5}}
\put(2,4){\circle*{.5}}
\put(6,4){\circle*{.5}}
\put(4,2){\circle*{.5}}
\put(8,2){\circle*{.5}}
\put(7,0){\circle*{.5}}
\put(9,0){\circle*{.5}}
\end{picture}
\,T14
\hfill \hfill
\begin{picture}(10,8)
\put(3,8){\line(-1,-1){2}}
\put(3,8){\line(1,-1){6}}
\put(5,6){\line(-1,-1){4}}
\put(3,4){\line(1,-2){1}}
\put(7,4){\line(-1,-2){1}}
\put(1,2){\line(-1,-2){1}}
\put(1,2){\line(1,-2){1}}
\put(9,2){\line(-1,-2){1}}
\put(9,2){\line(1,-2){1}}
\put(3,8){\circle*{.5}}
\put(1,6){\circle*{.5}}
\put(5,6){\circle*{.5}}
\put(3,4){\circle*{.5}}
\put(7,4){\circle*{.5}}
\put(1,2){\circle*{.5}}
\put(4,2){\circle*{.5}}
\put(6,2){\circle*{.5}}
\put(9,2){\circle*{.5}}
\put(0,0){\circle*{.5}}
\put(2,0){\circle*{.5}}
\put(8,0){\circle*{.5}}
\put(10,0){\circle*{.5}}
\end{picture}
\,T15
\hfill \hfill
\begin{picture}(8,8)
\put(1,8){\line(-1,-2){1}}
\put(1,8){\line(1,-1){6}}
\put(3,6){\line(-1,-1){2}}
\put(1,4){\line(-1,-2){1}}
\put(1,4){\line(1,-2){1}}
\put(5,4){\line(-1,-1){2}}
\put(3,2){\line(-1,-2){1}}
\put(3,2){\line(1,-2){1}}
\put(7,2){\line(-1,-2){1}}
\put(7,2){\line(1,-2){1}}
\put(1,8){\circle*{.5}}
\put(0,6){\circle*{.5}}
\put(3,6){\circle*{.5}}
\put(1,4){\circle*{.5}}
\put(5,4){\circle*{.5}}
\put(0,2){\circle*{.5}}
\put(2,2){\circle*{.5}}
\put(3,2){\circle*{.5}}
\put(7,2){\circle*{.5}}
\put(2,0){\circle*{.5}}
\put(4,0){\circle*{.5}}
\put(6,0){\circle*{.5}}
\put(8,0){\circle*{.5}}
\end{picture}
\,T16
\end{center}
\begin{center}
\setlength{\unitlength}{5pt}
\hfill
\begin{picture}(7,10)
\put(4,10){\line(-1,-1){2}}
\put(4,10){\line(1,-1){2}}
\put(6,8){\line(-1,-1){4}}
\put(6,8){\line(1,-2){1}}
\put(4,6){\line(1,-1){2}}
\put(2,4){\line(-1,-2){2}}
\put(2,4){\line(1,-2){1}}
\put(6,4){\line(-1,-2){1}}
\put(6,4){\line(1,-2){1}}
\put(1,2){\line(1,-2){1}}
\put(4,10){\circle*{.5}}
\put(2,8){\circle*{.5}}
\put(6,8){\circle*{.5}}
\put(4,6){\circle*{.5}}
\put(7,6){\circle*{.5}}
\put(2,4){\circle*{.5}}
\put(6,4){\circle*{.5}}
\put(1,2){\circle*{.5}}
\put(3,2){\circle*{.5}}
\put(5,2){\circle*{.5}}
\put(7,2){\circle*{.5}}
\put(0,0){\circle*{.5}}
\put(2,0){\circle*{.5}}
\end{picture}
\,T17
\hfill \hfill
\begin{picture}(7,10)
\put(4,10){\line(-1,-1){2}}
\put(4,10){\line(1,-1){2}}
\put(6,8){\line(-1,-1){4}}
\put(6,8){\line(1,-2){1}}
\put(4,6){\line(1,-1){2}}
\put(2,4){\line(-1,-2){1}}
\put(2,4){\line(1,-2){1}}
\put(6,4){\line(-1,-2){2}}
\put(6,4){\line(1,-2){1}}
\put(5,2){\line(1,-2){1}}
\put(4,10){\circle*{.5}}
\put(2,8){\circle*{.5}}
\put(6,8){\circle*{.5}}
\put(4,6){\circle*{.5}}
\put(7,6){\circle*{.5}}
\put(2,4){\circle*{.5}}
\put(6,4){\circle*{.5}}
\put(1,2){\circle*{.5}}
\put(3,2){\circle*{.5}}
\put(5,2){\circle*{.5}}
\put(7,2){\circle*{.5}}
\put(4,0){\circle*{.5}}
\put(6,0){\circle*{.5}}
\end{picture}
\,T18
\hfill \hfill
\begin{picture}(15,10)
\put(1,2){\line(1,-2){1}}
\put(1,2){\line(-1,-2){1}}
\put(5,2){\line(1,-2){1}}
\put(5,2){\line(-1,-2){1}}
\put(9,2){\line(1,-2){1}}
\put(9,2){\line(-1,-2){1}}
\put(13,2){\line(1,-2){1}}
\put(13,2){\line(-1,-2){1}}
\put(3,4){\line(1,-1){2}}
\put(3,4){\line(-1,-1){2}}
\put(11,4){\line(1,-1){2}}
\put(11,4){\line(-1,-1){2}}
\put(7,6){\line(2,-1){4}}
\put(7,6){\line(-2,-1){4}}
\put(14,6){\line(1,-2){1}}
\put(14,6){\line(-1,-2){1}}
\put(10,8){\line(-3,-2){3}}
\put(10,8){\line(2,-1){4}}
\put(7,10){\line(-3,-2){3}}
\put(7,10){\line(3,-2){3}}
\put(7,10){\circle*{.5}}
\put(4,8){\circle*{.5}}
\put(10,8){\circle*{.5}}
\put(7,6){\circle*{.5}}
\put(14,6){\circle*{.5}}
\put(3,4){\circle*{.5}}
\put(11,4){\circle*{.5}}
\put(13,4){\circle*{.5}}
\put(15,4){\circle*{.5}}
\put(1,2){\circle*{.5}}
\put(5,2){\circle*{.5}}
\put(9,2){\circle*{.5}}
\put(13,2){\circle*{.5}}
\put(0,0){\circle*{.5}}
\put(2,0){\circle*{.5}}
\put(4,0){\circle*{.5}}
\put(6,0){\circle*{.5}}
\put(8,0){\circle*{.5}}
\put(10,0){\circle*{.5}}
\put(12,0){\circle*{.5}}
\put(14,0){\circle*{.5}}
\end{picture}
\,T19
\hfill \makebox{}
\end{center}
\begin{center}
\setlength{\unitlength}{5pt}
\hfill
\begin{picture}(9,10)
\put(2,10){\line(-1,-1){2}}
\put(2,10){\line(1,-1){6}}
\put(4,8){\line(-1,-1){2}}
\put(6,6){\line(-1,-1){2}}
\put(4,4){\line(-1,-2){1}}
\put(4,4){\line(1,-2){2}}
\put(8,4){\line(-1,-2){1}}
\put(8,4){\line(1,-2){1}}
\put(5,2){\line(-1,-2){1}}
\put(2,10){\circle*{.5}}
\put(0,8){\circle*{.5}}
\put(4,8){\circle*{.5}}
\put(2,6){\circle*{.5}}
\put(6,6){\circle*{.5}}
\put(4,4){\circle*{.5}}
\put(8,4){\circle*{.5}}
\put(3,2){\circle*{.5}}
\put(5,2){\circle*{.5}}
\put(7,2){\circle*{.5}}
\put(9,2){\circle*{.5}}
\put(4,0){\circle*{.5}}
\put(6,0){\circle*{.5}}
\end{picture}
\,T20
\hfill \hfill
\begin{picture}(11,10)
\put(2,10){\line(-1,-1){2}}
\put(2,10){\line(1,-1){8}}
\put(4,8){\line(-1,-1){2}}
\put(6,6){\line(-1,-1){2}}
\put(8,4){\line(-1,-1){2}}
\put(4,4){\line(-1,-2){1}}
\put(4,4){\line(1,-2){1}}
\put(10,2){\line(-1,-2){1}}
\put(10,2){\line(1,-2){1}}
\put(2,10){\circle*{.5}}
\put(0,8){\circle*{.5}}
\put(4,8){\circle*{.5}}
\put(2,6){\circle*{.5}}
\put(6,6){\circle*{.5}}
\put(4,4){\circle*{.5}}
\put(8,4){\circle*{.5}}
\put(3,2){\circle*{.5}}
\put(5,2){\circle*{.5}}
\put(6,2){\circle*{.5}}
\put(10,2){\circle*{.5}}
\put(9,0){\circle*{.5}}
\put(11,0){\circle*{.5}}
\end{picture}
\,T21
\hfill \hfill
\begin{picture}(16,10)
\put(1,10){\line(-1,-2){1}}
\put(1,10){\line(2,-1){12}}
\put(5,8){\line(-2,-1){4}}
\put(1,6){\line(-1,-2){1}}
\put(1,6){\line(1,-2){1}}
\put(9,6){\line(-2,-1){4}}
\put(5,4){\line(-1,-1){2}}
\put(5,4){\line(1,-1){2}}
\put(13,4){\line(-1,-1){2}}
\put(13,4){\line(1,-1){2}}
\put(3,2){\line(-1,-2){1}}
\put(3,2){\line(1,-2){1}}
\put(7,2){\line(-1,-2){1}}
\put(7,2){\line(1,-2){1}}
\put(11,2){\line(-1,-2){1}}
\put(11,2){\line(1,-2){1}}
\put(15,2){\line(-1,-2){1}}
\put(15,2){\line(1,-2){1}}
\put(1,10){\circle*{.5}}
\put(0,8){\circle*{.5}}
\put(5,8){\circle*{.5}}
\put(1,6){\circle*{.5}}
\put(9,6){\circle*{.5}}
\put(0,4){\circle*{.5}}
\put(2,4){\circle*{.5}}
\put(5,4){\circle*{.5}}
\put(13,4){\circle*{.5}}
\put(3,2){\circle*{.5}}
\put(7,2){\circle*{.5}}
\put(11,2){\circle*{.5}}
\put(15,2){\circle*{.5}}
\put(2,0){\circle*{.5}}
\put(4,0){\circle*{.5}}
\put(6,0){\circle*{.5}}
\put(8,0){\circle*{.5}}
\put(10,0){\circle*{.5}}
\put(12,0){\circle*{.5}}
\put(14,0){\circle*{.5}}
\put(16,0){\circle*{.5}}
\end{picture}
\,T22
\hfill \makebox{}
\end{center}

\end{appendix}

\bibliography{quantum}
\bibliographystyle{mybib}

\end{document}